\documentclass[12pt]{article}
\usepackage{amssymb,amsmath,amsthm,amsfonts,mathrsfs}
\usepackage{subfigure,epsfig,psfrag,epstopdf}
\usepackage{fullpage,graphicx,color}
\usepackage{algorithm,algorithmic}
\usepackage{tabularx,array,booktabs,bm,multirow}
\usepackage{longtable}
\usepackage{cite,url}
\usepackage[bookmarks=false]{hyperref}
\usepackage[labelfont=bf]{caption}

\newcommand{\mB}{\mathcal{B}}
\newcommand{\mD}{\mathcal{D}}
\newcommand{\mJ}{\mathcal{J}}

\newcommand{\mC}{\mathcal{C}}

\newcommand{\mP}{\mathcal{P}}
\newcommand{\RR}{\mathbb{R}}

\def\bnabla{\boldsymbol{\nabla}}

\def\Dpartial#1#2{\dfrac{\partial #1}{\partial #2}}
\newcommand{\argmin}{\operatorname{argmin}}

\newcommand*{\mybox}[1]{\framebox{#1}}

\begin{document}

\title{On Reconstruction of Binary Images by Efficient Sample-based Parameterization in
Applications for Electrical Impedance Tomography}

\date{}

\author{{\bf \normalsize Paul R.~Arbic II} \\ {\it \small Department of Mathematical Sciences,
Florida Institute of Technology, Melbourne, FL 32901, USA} \and
{\bf \normalsize Vladislav Bukshtynov}\footnote{Corresponding author: \texttt{vbukshtynov@fit.edu}} \\
{\it \small Department of Mathematical Sciences, Florida Institute of Technology, Melbourne, FL 32901, USA}}
\maketitle

\begin{abstract}
An efficient computational approach for optimal reconstruction of binary-type images
suitable for models in various applications including biomedical imaging is
developed and validated. The methodology includes derivative-free
optimization supported by a set of sample solutions with customized
geometry generated synthetically. The reduced dimensional control space
is organized based on contributions from individual samples and the
efficient parameterization obtained from the description of the samples'
geometry. The entire framework has an easy-to-follow design due to
a nominal number of tuning parameters which makes the approach simple
for practical implementation in various settings, as well as for
adjusting it to new models and enhancing the performance. High efficiency
in computational time is achieved through applying the coordinate descent
method to work with individual controls in the predefined custom order.
This technique is shown to outperform commonly used gradient-based and other
derivative-free methods with applied PCA-based control space reduction in terms of
both qualities of binary images and stability of obtained solutions when noise is
present in the measurement data. Performance of the complete computational
framework is tested in applications to 2D inverse problems of detecting
inclusions or defective regions by the electrical impedance tomography~(EIT).
The results demonstrate the efficient performance of the new method and its high
potential for improving the overall quality of the EIT-based procedures
including biomedical applications.

\end{abstract}

\section{Introduction}
\label{sec:introduction}

In this work, we propose a novel computational approach for optimal
reconstruction of biomedical images based on any available measurements
usually obtained with some noise. In particular, this approach is useful
in various applications for medical practices dealing with models
characterized by parameters with or near to binary-type distributions,
e.g., heat or electrical conductivity. The proposed computational framework
is built around a derivative-free optimization algorithm and supported
by a set of sample solutions. These samples are generated synthetically
with a geometry based on any available prior knowledge of the simulated
phenomena and the expected structure of obtained images. The ease of
parallelization allows operations on very large sample sets which
enables the best approximations for the initial guess and, as a result,
close proximity to the local/global solution of the optimal control
problem. The controls are effectively defined to utilize individual
contributions from samples selected to represent the compound image
and the efficient parameterization obtained from the description of
the samples' geometry. As our new approach shows similarity to
many modern learning algorithms, it could be included in the category
of dictionary learning where the computed solutions are used as a training
set for the learning process, see \cite{Fan2020,Hamilton2018,Wei2019,Ren2020}
for more details and parallels in this area.

The proposed computational framework has an easy-to-follow design and
is tuned by a nominal number of computational parameters making the
approach simple for practical implementation for various applications
also beyond biomedical imaging. High computational efficiency is
achieved by applying the coordinate descent method customized to work
with individual controls in the predefined custom order.

As known from practical applications, fine scale optimization performed
on fine meshes provides high resolution images. On the other hand,
such solutions will require computational time increased due to the
size of the fine mesh and the associated control space. In addition,
obtained images may not provide clear boundaries between regions
identified by different physical properties in space. As a result,
a smooth transition cannot provide an accurate recognition of shapes,
e.g., of cancer-affected regions while solving an inverse problem of
cancer detection~(IPCD). In our computations, fine mesh is used only
to assess the measurement fit in terms of evaluated cost functionals.

The initial motivation for the research described in the current
paper was to apply our new computational approach to IPCD by the
electrical impedance tomography~(EIT) technique. However, this
methodology could be easily extended to a broad range of problems
in biomedical sciences, also in physics, geology, chemistry, etc.
EIT is a rapidly developing non-invasive imaging technique gaining
popularity by enabling various medical applications to perform
screening for cancer detection \cite{Brown2003,Holder2004,Lionheart2004,
Abascal2011,Choi2007,Boverman2008,Adler2008,Uhlmann2009,Zou2003Review,Liu2019}.
A well-known fact is that the electrical properties, e.g., electrical
conductivity or permittivity, of different tissues are different
if they are healthy or affected by cancer. This phenomenon is used
in EIT to produce images of biological tissues by interpreting their
response to applied voltages or injected currents
\cite{Brown2003,Holder2004,Cheney1999}. The inverse EIT problem
deals with reconstructing the electrical conductivity by measuring
voltages or currents at electrodes placed on the surface of a test
volume. This so-called Calderon-type inverse problem \cite{Calderon1980}
is highly ill-posed, refer to topical review \cite{Borcea2002}.
Since the 1980s, various techniques have been suggested to solve it
computationally. We refer to the recent papers
\cite{Adler2015,Bera2018,Wang2020} with review on the current state
of the art and the existing open problems associated with EIT and
its applications. Generally speaking, the EIT-based methodology
could be easily applied in various fields where nondestructive testing
is needed to identify ``defective regions'' or regions with properties
considered abnormal in comparison with expected (normal) values associated
with the rest of the domain. To take advantage of this flexibility, in
the rest of the paper, we refer to the defects, also inclusions or targets,
in the material (media) to be identified using the proposed approach within
the general framework of EIT applications.

This paper proceeds as follows. In Section~\ref{sec:math}, we present
a very general mathematical description of the inverse EIT problem
formulated as an optimal control problem. Procedures for solving this
optimization problem with the proposed sample-based parameterization
are discussed in Section~\ref{sec:solution}. Model descriptions and
detailed computational results including discussion on chosen methods
are presented in Section~\ref{sec:results}. Concluding remarks are
provided in Section~\ref{sec:remarks}.

\section{Mathematical Model for Inverse EIT Problem}
\label{sec:math}

As discussed at length in \cite{AbdullaBukshtynovSeif2021} and \cite{KoolmanBukshtynov2021},
the inverse EIT problem is formulated as a partial differential equation (PDE) constrained
optimal control problem for an open and bounded set $\Omega \subset \RR^n, \ n = 2, 3$,
representing body with electrical conductivity at point $x \in \Omega$ given by function
$\sigma(x): \, \Omega \rightarrow \RR_+$. In this paper, we use
the so-called ``voltage--to--current" model where voltages (electrical
potentials) $U = (U_{\ell})_{\ell=1}^m \in \RR^m$ are applied to $m$
electrodes $(E_{\ell})_{\ell=1}^m$ with contact impedances
$(Z_{\ell})_{\ell=1}^m \in \RR^m_+$ subject to the ground
(zero potential) condition
\begin{equation}
  \sum_{\ell=1}^m U_{\ell}=0.
  \label{eq:ground_conds}
\end{equation}
These voltages initiate electrical currents
$(I_{\ell})_{\ell=1}^m\in \RR^m$ through the same electrodes $E_{\ell}$
placed at the periphery of the body $\partial \Omega$.
The electrical currents may be computed as
\begin{equation}
  I_{\ell} = \int_{E_{\ell}}  \sigma(x) \Dpartial{u(x)}{n} \, ds,
  \quad \ell = 1, \ldots, m
  \label{eq:el_current}
\end{equation}
based on conductivity field $\sigma(x)$ and a distribution of electrical
potential $u(x): \, \Omega \rightarrow \RR$ obtained as a solution of the
following elliptic problem:
\begin{subequations}
  \begin{alignat}{3}
    \bnabla \cdot \left[ \sigma(x) \bnabla u(x) \right] &= 0,
    && \quad x \in \Omega \label{eq:forward_1}\\
    \Dpartial{u(x)}{n} &= 0, && \quad x \in \partial \Omega -
    \bigcup\limits_{\ell=1}^{m} E_{\ell} \label{eq:forward_2}\\
    u(x) +  Z_{\ell} \sigma(x) \Dpartial{u(x)}{n} &= U_{\ell},
    && \quad x \in E_{\ell}, \ \ell= 1, \ldots, m
    \label{eq:forward_3}
  \end{alignat}
  \label{eq:forward}
\end{subequations}
in which $n$ is an external unit normal vector on $\partial \Omega$.
The complete description and analysis of electrode models used in
electric current computed tomography may be found in \cite{Somersalo1992}.

We set conductivity $\sigma(x)$ here as a control variable and
formulate the inverse EIT (conductivity) problem \cite{Calderon1980}
as a PDE-constrained optimal control problem \cite{AbdullaBukshtynovSeif2021}
by considering least-square minimization of mismatches
$\left( I_{\ell} - I_{\ell}^* \right)^2$, where
$(I_{\ell}^*)_{\ell=1}^m \in \RR^m$ are measurements of electrical
currents $I_{\ell}$. In addition, we have to mention a well-known
fact that this inverse EIT problem to be solved in a discretized
domain $\Omega$ is highly ill-posed. Therefore, we enlarge the data
up to size of $m^2$ by adding new measurements following the
``rotation scheme'' described in detail in \cite{AbdullaBukshtynovSeif2021}
while keeping the size of the unknown parameters, i.e., elements in the
discretized description for $\sigma(x)$, fixed. Having a new set of
data $(I_{\ell}^{k*})_{\ell,k=1}^m$ and in light of the Robin
condition \eqref{eq:forward_3} used together with \eqref{eq:el_current},
we define a complete form of the cost functional
\begin{equation}
  \mJ (\sigma) = \sum_{k=1}^{m} \sum_{\ell=1}^m
  \left[ \int_{E_{\ell}} \dfrac{U^k_{\ell}-u^k(x;\sigma)}{Z_{\ell}}
  \, ds - I^{k*}_{\ell} \right]^2
  \label{eq:cost_functional}
\end{equation}
for the optimal control problem
\begin{equation}
  \hat \sigma(x) = \underset{\sigma}{\argmin} \ \mJ(\sigma)
  \label{eq:minJ_sigma}
\end{equation}
subject to PDE constraint \eqref{eq:forward} where each function
$u^k(\cdot; \sigma), \ k = 1, \ldots, m$, solves elliptic PDE problem
\eqref{eq:forward_1}--\eqref{eq:forward_3}. We also note that these
solutions of the forward EIT problem after applying \eqref{eq:el_current}
and adding some noise may be used for generating various model examples
(synthetic data) for inverse EIT problems to adequately mimic the
presence of defective regions mentioned above.

\section{Solution by Sample-based Parameterization}
\label{sec:solution}

\subsection{Preliminaries and Main Notations}
\label{sec:step_0}

Without loss of generality, here we discuss our new algorithm for
solving problem \eqref{eq:minJ_sigma} in 2D ($n = 2$) domain
\begin{equation}
  \Omega = \left\{ x \in \RR^2 : \ | x |^2 < R^2 \right\}
  \label{eq:domain}
\end{equation}
which is a disc of radius $R$. However, the same analysis could be
easily extended to any complexity of domain $\Omega$ in 3D ($n = 3$)
settings. In addition, we assume that the actual (true) electrical
conductivity $\sigma_{true}(x)$ we seek to reconstruct could be
represented by
\begin{equation}
  \sigma_{true}(x) = \left\{
  \begin{aligned}
    \sigma_d, & \quad x \in \Omega_d,\\
    \sigma_n, & \quad x \in \Omega_n,
  \end{aligned}
  \right.
  \quad \Omega_d \cap \Omega_n = \emptyset
  \label{eq:sigma_true}
\end{equation}
where $\sigma_d$ and $\sigma_n$ are known constants for the respective
defective region $\Omega_d$ and the rest of domain $\Omega_n$ with normal
(expected) conductivity.

\begin{figure}[htb]
  \begin{center}
  \includegraphics[width=0.33\textwidth]{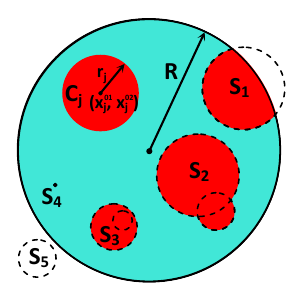}
  \end{center}
  \caption{Different scenarios of $j$th circle appearance in the
    $i$th sample: regular case $C_j$ and special cases $S_1$
    through $S_5$.}
  \label{fig:sample_param}
\end{figure}

We seek for the solution of \eqref{eq:minJ_sigma} in a form
\begin{equation}
  \sigma(x) = \sum_{i=1}^{N_s} \alpha_i \bar \sigma_i(x),
  \qquad 0 \leq \alpha_i \leq 1,
  \label{eq:sigma_main}
\end{equation}
where $\bar \sigma_i(x)$, $i = 1, \ldots, N_s$, are sample solutions
generated synthetically and convexly weighted by coefficients
$\alpha_i$
\begin{equation}
  \sum_{i=1}^{N_s} \alpha_i = 1.
  \label{eq:sigma_alpha}
\end{equation}
The entire collection of $N$ samples
\begin{equation}
  \mC(N) = (\bar \sigma_i(x))_{i=1}^{N}, \qquad N \gg N_s
  \label{eq:smpl_collect}
\end{equation}
in fact could be generated based on any assumptions made for the
(geometrical) structure of the reconstructed images. Here we
assume that clear shapes for binary images could be obtained by
combining simple convex geometric shapes (elements) in 2D such
as triangles, squares, circles, etc. For example, in the current
research the $i$th sample in our $N$-collection consists of $N_c^i$
circles of various radii $r \in \RR_+$ and centers
$x^0 = (x^{01}, x^{02}) \in \RR^2$ located inside domain $\Omega$, i.e.
\begin{equation}
  \bar \sigma_i(x) = \left\{
  \begin{aligned}
    \sigma_d, &\quad | x - x^0_j  |^2 \leq r^2_j, \ j = 1,
    \ldots, N_c^i\\
    \sigma_n, &\quad {\rm otherwise}
  \end{aligned}
  \right.
  \label{eq:smpl_par}
\end{equation}
Thus, some approximations of $\sigma_d$ and $\sigma_n$ in \eqref{eq:sigma_true}
to be used then in \eqref{eq:smpl_par} are required and considered as {\it a priori} knowledge
needed to apply the approach in practice. In \eqref{eq:smpl_par}, all $N_c^i$ circles
parameterized by the set of triplets
\begin{equation}
  \mP_i = (\{ x^{01}_j, x^{02}_j, r_j\})_{j=1}^{N_c^i},
  \qquad i = 1, \ldots, N
  \label{eq:smpl_triplet}
\end{equation}
are generated randomly subject to the following restrictions:
\begin{subequations}
  \begin{alignat}{3}
    &&|x^0_j| < R+r_j, \qquad j &= 1, \ldots, N_c^i,
    \label{eq:smpl_restr_1}\\
    &&1 \leq N_c^i \leq N_{c,\max}, \qquad i &= 1, \ldots, N.
    \label{eq:smpl_restr_2}
  \end{alignat}
  \label{eq:smpl_restr}
\end{subequations}
Parameter $N_{c,\max}$ in \eqref{eq:smpl_restr_2} defines the
maximum number of circles in the samples and, in fact, sets the
highest level of complexity (resolution) for the reconstructed
image $\hat \sigma(x)$. Figure~\ref{fig:sample_param} shows
different scenarios of the $j$th circle's appearance in the
$i$th sample: regular case $C_j$ (fully inside $\Omega$) and
a few special cases.
\begin{itemize}
  \item[(a)] $S_1$ for circles which are partially outside the
    domain $\Omega$;
  \item[(b)] $S_2$ and $S_3$ for circles with respective partially
    and fully overlapped regions; and
  \item[(c)] degenerate cases $S_4$ and $S_5$ correspondingly
    of zero radius or appear fully outside of domain $\Omega$.
\end{itemize}
We note that all circles of the special cases mentioned in (c) are
rejected when the samples are generated.

After completing the collection of $N$ sample solutions following
the description above, the proposed computational algorithm for
solving problem \eqref{eq:minJ_sigma} could be executed in two steps:
\begin{itemize}
  \item[{\bf Step~1:}] Define the initial basis of samples
    \begin{equation}
      \mB^0 = (\bar \sigma_i(x))_{i=1}^{N_s}
      \label{eq:smpl_basis}
    \end{equation}
    by choosing $N_s$ best samples out of collection $\mC(N)$ which
    provide the best measurement fit in terms of cost functional \eqref{eq:cost_functional}.
  \item[{\bf Step~2:}] Set all parameters $(\mP_i)_{i=1}^{N_s}$ in the
    description of sample basis $\mB$ as controls to perform
    optimization for solving problem \eqref{eq:minJ_sigma} numerically
     to find the optimal basis $\hat \mB$.
\end{itemize}

\subsection{Step~1: Defining Initial Basis}
\label{sec:step_1}

This step requires solving forward problem \eqref{eq:forward} and
evaluating cost functional \eqref{eq:cost_functional} $N$ times for
all samples in the $N$-collection. For a fixed scheme of potentials
$U$ the data $\mD_i = I^{k}_{\ell}(\sigma_i) \in \RR^{m^2}, \ i = 1,
\ldots, N$, could be precomputed by \eqref{eq:forward} and
\eqref{eq:el_current} and then stored for multiple uses with different
models. In addition, this task may be performed in parallel with minimal
computational overhead which allows easy switching between various
schemes for electrical potentials. Easy parallelization enables taking
$N$ quite large which helps better approximate the solution by the
initial state of basis $\mB$ before proceeding to Step~2.

The number of samples $N_s$ in basis $\mB$ may be defined experimentally
based on the model complexity. We suggest $N_s$ to be sufficiently
large to properly support a local/global search for optimal solution
$\hat \sigma(x)$ during Step~2. At the same time this number should
allow the total number of controls while solving problem
\eqref{eq:minJ_sigma} to be comparable with data dimension, namely
$m^2$, for satisfying the well-posedness requirement.

Working with the models of complicated structures may require
increasing the current number of elements (circles) in every sample
within the chosen basis $\mB$. In such a case, one could re-set
parameters $N_c^i, \ i = 1, \ldots, N_s$, to higher values and add
missing elements, for example, by generating new circles randomly
as degenerate cases $S_4$. This will project the initial basis
$\mB^0$ onto a new control space of a higher dimension without
any loses in the quality of the initial solution.

Step~1 will be completed after ranking samples in ascending order in
terms of computed cost functionals \eqref{eq:cost_functional} while
comparing the obtained data $(\mD_i)_{i=1}^{N}$ with true data
$(I_{\ell}^{k*})_{\ell,k=1}^m$ available from the actual measurements.
After ranking, the first $N_s$ samples will create the initial basis
$\mB^0$ to be used in place of the initial guess for optimization
in Step~2.

\subsection{Step~2: Solving Optimal Control Problem}
\label{sec:step_2}

As discussed in Section~\ref{sec:step_0}, all elements (circles) in
all samples of basis $\mB$ obtained during Step~1 ranking procedure
will be represented by a finite number of ``sample-based'' parameters
$(\mP_i)_{i=1}^{N_s}$. In general, solution $\sigma(x)$ could be
uniquely represented as a function of $\mP_i, \ i = 1, \ldots, N_s$.
The continuous form of optimal control problem \eqref{eq:minJ_sigma}
may be substituted with its new equivalent form defined over the
finite set of controls $\mP = (\mP_i)_{i=1}^{N_s}$. In addition to
this, problem \eqref{eq:minJ_sigma} could be further extended by
adding weights $\alpha = (\alpha_i)_{i=1}^{N_s}$
in \eqref{eq:sigma_main} to the set of new controls. After this we
arrive at the final form of the optimization problem to be later
solved numerically:
\begin{equation}
  (\hat \mP, \hat \alpha) = \underset{\mP,\alpha}{\argmin}
  \ \mJ(\mP,\alpha)
  \label{eq:minJ_ext}
\end{equation}
subject to PDE constraint \eqref{eq:forward}, linear constraint
\eqref{eq:sigma_alpha} and properly established bounds for all
components of control $(\mP, \alpha)$.
As easily followed from the structure of this new control, a
dimension of the parameterized solution space is bounded by
\begin{equation}
  \max(\dim(\mP, \alpha)) = N_s \cdot [ N_{c,\max} (n+1) + 1].
\end{equation}

When solving \eqref{eq:minJ_ext} iteratively, one may choose
to terminate the optimization run at $k$th (major) iteration
once the following criterion is satisfied:
\begin{equation}
    \left| \dfrac{\mJ^k - \mJ^{k-1}}{\mJ^k}
    \right| < \epsilon
  \label{eq:termination}
\end{equation}
subject to chosen tolerance $\epsilon \in \RR_+$. Although both
\eqref{eq:minJ_sigma} and \eqref{eq:minJ_ext} are obviously
not separable optimization problems, the coordinate descent (CD)
method is used to solve \eqref{eq:minJ_ext}. This choice is
motivated by several reasons, namely
\begin{itemize}
  \item simplicity of the form for establishing the equivalence
    between controls $\sigma(x)$ and $(\mP, \alpha)$ provided by
    \eqref{eq:sigma_main} and
    \eqref{eq:smpl_par}--\eqref{eq:smpl_triplet},
  \item close proximity of samples in the initial basis $\mB^0$
    to the local/global solutions after completing Step~1, and
  \item straightforward computational implementation.
\end{itemize}
The efficiency of the entire optimization framework is confirmed
by extensive computational results for multiple models of different
complexity presented in Section~\ref{sec:results}. A summary of the
complete computational framework to perform our new optimization with
sample-based parameterization is provided in Algorithm~\ref{alg:main_opt}.
We also note that in order to improve the computational efficiency,
the applied CD method is modified by specifying the order in which all
controls are perturbed while solving problem \eqref{eq:minJ_ext} as
discussed in detail in Section~\ref{sec:CD_method}.
Briefly, instead of following the sensitivity feedback, choosing
controls is considered in the sample-by-sample order, and within each
sample $\bar \sigma_i$ we optimize over all circles' triplets
$\{ x^{01}_j, x^{02}_j, r_j\}$ and then over the sample's weight
$\alpha_i$, see Step~2 in Algorithm~\ref{alg:main_opt} for clarity.

\begin{algorithm}[htb!]
\begin{algorithmic}
  \STATE set parameters: $N$, $N_{c,\max}$, $N_s$
  \FOR{$i \leftarrow 1$ to $N$}
    \STATE generate $\bar \sigma_i(x)$ by
    \eqref{eq:smpl_par}--\eqref{eq:smpl_restr}
    \STATE obtain data $\mD_i = I^{k}_{\ell}(\bar \sigma_i)$ by
      \eqref{eq:forward} and \eqref{eq:el_current}
  \ENDFOR
  \STATE \begin{center} \mybox{Step~1}
         \it{Defining Initial Basis} \end{center}
  \STATE select model and obtain true data
    $(I_{\ell}^{k*})_{\ell,k=1}^m$
  \FOR{$i \leftarrow 1$ to $N$}
    \STATE compute $\mJ(\bar \sigma_i)$ by \eqref{eq:cost_functional}
  \ENDFOR
  \STATE choose $N_s$ best samples from $\mC(N)$ by values
    $\mJ(\bar \sigma_i)$
  \STATE form initial basis $\mB^0$ subject to $N_{c,\max}$
  \STATE set initial weights $\alpha^0$
  \STATE compute $\sigma^0(x)$ using $\mB^0$ and $\alpha^0$ by
    \eqref{eq:sigma_main}
  \STATE \begin{center} \mybox{Step~2}
         \it{Solving Optimal Control Problem} \end{center}
  \STATE $k \leftarrow 0$
  \REPEAT
  \FOR{$i \leftarrow 1$ to $N_s$}
    \FOR{$j \leftarrow 1$ to $N^{i}_c$}
      \STATE {\bf optimize} over $j$th circle triplet
        $\{ x^{01}_j, x^{02}_j, r_j\}$ in $\bar \sigma_i$
    \ENDFOR
    \STATE {\bf optimize} over weight $\alpha_i$
  \ENDFOR
  \STATE $k \leftarrow k + 1$
  \STATE update $\sigma^k(x)$ using new basis $\mB^k$ and
    weights $\alpha^k$ by \eqref{eq:sigma_main}
  \UNTIL termination criterion \eqref{eq:termination} is satisfied
    to given tolerance
\end{algorithmic}
\caption{Computational workflow for optimization with sample-based
         parameterized controls}
\label{alg:main_opt}
\end{algorithm}

\subsection{Customized Coordinate Descent Method}
\label{sec:CD_method}

To conclude on the practical implementation of our modified
coordinate descent approach, customarily adopted for the
proposed optimization framework, here we discuss some details.
We also note that, despite the superior performance of the CD
approach, its current version used for all computations in this
paper still has enough room for further improvements towards
computational efficiency and applicability for 2D and 3D problems.

After finishing Step~1 (refer to Section~\ref{sec:step_1} for details),
we then have a multitude of potential methods for performing optimization
at Step~2 (see Section~\ref{sec:step_2}). As mentioned before, our
preferred method is a customized version of the CD algorithm. Coordinate
descent, or alternating variables, algorithms are a form of
discrete gradient algorithms that determine if there is an apparent
gradient by comparing the difference between the values of the
cost functional when measured at two or more positions in the related domains,
see \cite{Nocedal2006} for more details. Gradient descent is then simulated
by moving in discrete intervals, rather than a distance proportionate to
the estimated gradient's size as is common with typical gradient descent
techniques.

In our algorithm, we start with initializing controls obtained from
$N_s$ samples in use, namely all parameters $(\mP_i)_{i=1}^{N_s}$
from the initial basis $\mB^0 = (\bar \sigma^0_i(x))_{i=1}^{N_s}$
and associated weights $(\alpha^0_i)_{i=1}^{N_s}$. We note that
$\mB^0$ and $\alpha^0$ are used to compute $\sigma^0(x)$, and the
composite sample (image) to represent the final solution for Step~1
and simultaneously the initial guess used in Step~2.

The optimization process by the CD approach at Step~2 chooses
controls one-by-one from the entire control space
\begin{equation}
  (\mP,\alpha) = \{(\mP_i)_{i=1}^{N_s} \cup (\alpha _i)_{i=1}^{N_s}\}
  \label{eq:ctrl_space}
\end{equation}
in the sample--by--sample order and circle-by-circle order
within a sample. The ordered circles in all samples are sequentially
moved horizontally, then vertically and finally resized by changing
their radii. After completing the optimization over all circles in
the $i$th sample, its weight $\alpha_i$ is also optimized to employ
the benefits of the updated geometry of the current sample.

Now we focus on a single (minor)
optimization iteration within this procedure for changing
a single control $\zeta$ selected from $(\mP,\alpha)$. A summary
for all steps within this iteration is provided in
Algorithm~\ref{alg:cd_opt}.
\begin{algorithm}[htb!]
\begin{algorithmic}
  \STATE choose control $\zeta \in \{(\mP_i)_{i=1}^{N_s} \cup
    (\alpha _i)_{i=1}^{N_s}\}$
  \STATE compute $\mJ_{curr}(\bar \sigma_i)$
    by \eqref{eq:cost_functional}
  \STATE $d \leftarrow 1$
  \STATE $s \leftarrow 0$
  \REPEAT
    \IF{$d > 0$}
      \STATE increment $\zeta$ by \eqref{eq:smpl_triplet2}
        or \eqref{eq:smpl_triplet3}--\eqref{eq:smpl_triplet4}
    \ELSE
      \STATE decrement $\zeta$ by \eqref{eq:smpl_triplet2}
        or \eqref{eq:smpl_triplet3}--\eqref{eq:smpl_triplet4}
    \ENDIF
    \STATE compute $\mJ_{temp}(\bar \sigma_i)$ by \eqref{eq:cost_functional}
    \IF{$\mJ_{curr} \leq \mJ_{temp}$}
      \IF{s = 0}
        \STATE $d \leftarrow -1$
        \STATE decrement $\zeta$ by \eqref{eq:smpl_triplet2}
          or \eqref{eq:smpl_triplet3}--\eqref{eq:smpl_triplet4}
        \STATE $\mJ_{temp} \leftarrow \mJ_{curr}$
      \ENDIF
      \STATE $s \leftarrow s + 1$
      \IF{$s = s_{\max}$}
        \STATE \bf{break}
      \ENDIF
    \ELSE
      \STATE $s \leftarrow 1$
    \ENDIF
  \UNTIL{$\mJ_{curr} < \mJ_{temp}$}
\end{algorithmic}
\caption{CD method: minor iteration}
\label{alg:cd_opt}
\end{algorithm}

First, we choose a single parameter $\zeta$, for example, from the
set of controls $\mP$, as the current control variable and perturb
it by a preset value, say constant step-size~$\gamma$:
\begin{equation}
  \zeta^{new} = \zeta^{old} \pm \gamma, \qquad
  \zeta \in (\mP_i)_{i=1}^{N_s}.
  \label{eq:smpl_triplet2}
\end{equation}
Then we determine if the cost functional $\mJ(\bar \sigma)$ of the
altered composite sample $\bar \sigma(x)$ improves (decreases) or
worsens (increases). We always perturb the control $\zeta$ in the
same direction along the given axis initially and then, depending
upon if the cost functional improves or not, we decide how to proceed.
If the cost functional improves, we perturb it in the same direction again
and continue in that direction until we see an increase in the
functional value. At this point, we select the best position found
and move to the next control. If the initial perturbation worsens
the functional, we then reverse direction along the given axis and
continue in the same manner as before until we reach a point at which
the cost functional worsens.

Perturbations in both directions may result in a position/radius
associated with the cost functional value which is worse than its
value computed at the initial position/radius. In this case,
we revert to the parameter's initial value and continue to the
next parameter (control variable). If there is no change in the
functional value, we continue in the same direction until there
is either a change in the functional or some pre-specified amount
of steps $s_{max}$ is reached, at which point the same behavior as
described before is continued or we move to the first position that
had the same functional value, respectively.

We perform a highly similar process when optimizing over the weights
of each sample image in the composite, $\zeta \in \alpha$, but with
a key difference. Being that instead of perturbing $\zeta = \alpha_i$
by some $\gamma$ value, we multiply the weight $\alpha_i$ by
$1 \pm \gamma$, i.e.
\begin{equation}
  \zeta^{new} = (1 \pm \gamma) \, \zeta^{old}, \qquad
  \zeta \in (\alpha_i)_{i=1}^{N_s},
  \label{eq:smpl_triplet3}
\end{equation}
and then normalize all weights subject to the convexity condition
\eqref{eq:sigma_alpha}:
\begin{equation}
  \alpha_{i}^{new} =
  \dfrac{\alpha_{i}^{old}}{1 \pm \gamma \zeta^{old}},
  \qquad i = 1, \ldots, N_s.
  \label{eq:smpl_triplet4}
\end{equation}

We would like to conclude on our modified CD approach by re-iterating
that all operations across the controls are performed following the
specific ``schedule''. We start with the first circle in the description
of the best matching sample identified during the Step~1 procedure.
Then we proceed towards each consecutive control in this sample, namely
$x^{01}_{1j}$ followed by $x^{02}_{1j}$ followed by $r_{1j}$,
$j = 1, \ldots, N_c^i$, until every parameter is optimized at which
point the weight $\alpha_1$ for that sample is optimized. Then we proceed
to the next sample, $i = 2, 3, \ldots, N_s$, and do the same. Once all
$N_s$ samples have been optimized, we loop back to the first sample image
and perform the same operations again and continue this process until
at least one of the termination criteria is met.

\section{Computational Results}
\label{sec:results}

\subsection{Computational Model in 2D}
\label{sec:comp_model}

Our optimization framework integrates computational facilities for
solving forward PDE problem \eqref{eq:forward} and evaluating
cost functionals by \eqref{eq:cost_functional}. These facilities
are incorporated by using {\tt FreeFem++}, see \cite{FreeFem2012}
for details, an open-source, high-level integrated development
environment for obtaining numerical solutions of PDEs based on the
finite element nethod (FEM). For solving numerically forward PDE
problem \eqref{eq:forward}, spatial discretization is carried out
by implementing 7730 FEM triangular finite elements: P2 piecewise
quadratic (continuous) representation for electrical potential $u(x)$
and P0 piecewise constant representation for conductivity field
$\sigma(x)$. Systems of algebraic equations obtained after such
discretization are solved with {\tt UMFPACK}, a solver for
nonsymmetric sparse linear systems \cite{UMFPACK}.

\begin{figure}[htb!]
  \begin{center}
  \mbox{
  \subfigure[]{\includegraphics[width=0.5\textwidth]{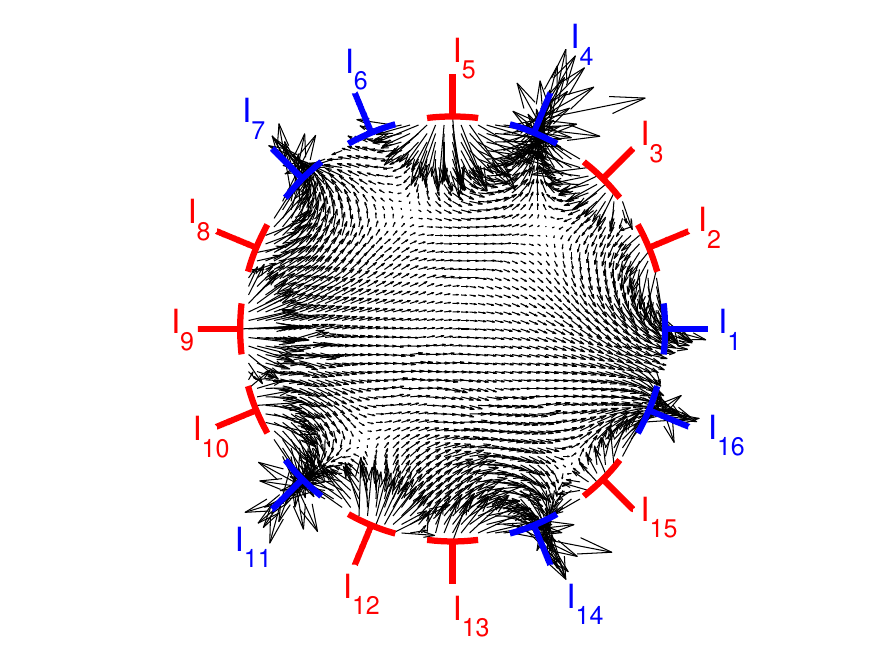}}
  \subfigure[]{\includegraphics[width=0.5\textwidth]{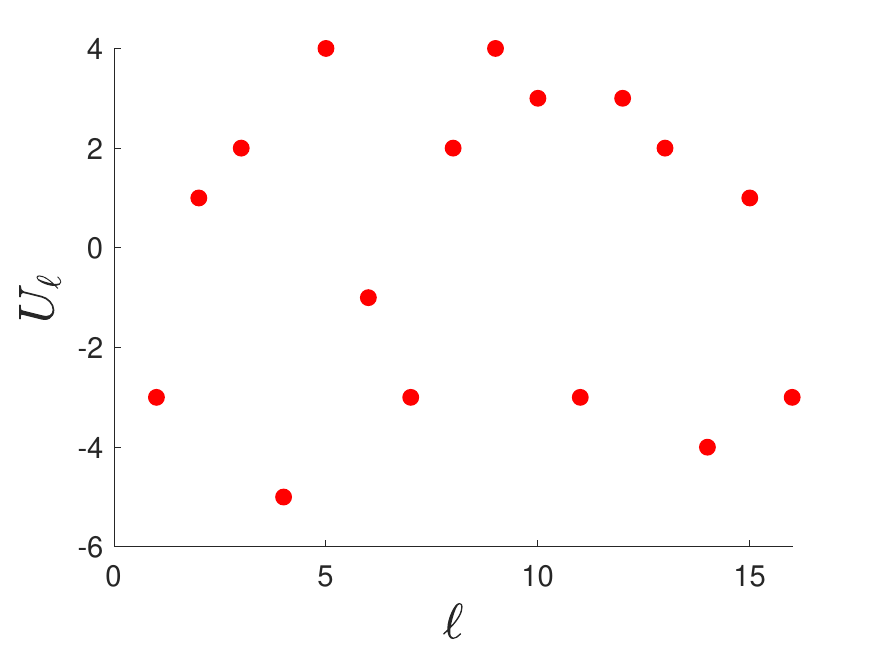}}}
  \end{center}
  \caption{(a)~Equispaced geometry of electrodes $(E_{\ell})_{\ell=1}^{16}$
    and electrical currents $I_{\ell}$ (positive in red, negative in blue)
    measured at $E_{\ell}$. Black arrows show the distribution of flux
    $\sigma(x) \bnabla u(x)$ in the interior of $\Omega$.
    (b)~Electrical potentials $U_{\ell}$.}
  \label{fig:model}
\end{figure}

All computations are performed using a 2D domain \eqref{eq:domain}
which is a disc of radius $R = 0.1$ with $m = 16$ equidistant electrodes
$E_{\ell}$ with half-width $w = 0.12$ rad covering approximately 61\%
of boundary $\partial \Omega$ as shown in Figure~\ref{fig:model}(a).
Electrical potentials $U_{\ell}$, see Figure~\ref{fig:model}(b), are
applied to electrodes $E_{\ell}$ following the ``rotation scheme''
discussed in Section~\ref{sec:math} and chosen to be consistent with
the ground potential condition \eqref{eq:ground_conds}. Determining
the Robin part of the boundary conditions in \eqref{eq:forward_3}, we
equally set the electrode contact impedance $Z_{\ell} = 0.1$.
To solve optimal control problem \eqref{eq:minJ_ext} iteratively, our
framework is utilizing derivative-free Coordinate Descent or
alternating variables approach, see \cite{Nocedal2006} for more details.
The actual (true) electrical conductivity $\sigma_{true}(x)$ we seek to
reconstruct is defined analytically for each model in \eqref{eq:sigma_true}
by setting $\sigma_d = 0.4$ and $\sigma_n = 0.2$. The initial guess for
control $\mP$ at Step~2 is provided by the parameterization of initial
basis $\mB^0$ obtained after completion of Step~1. For control $\alpha$,
the initial values are set to be equal, i.e., $\alpha^0_i = 1/N_s$.
Termination criteria are set by tolerance $\epsilon = 10^{-4}$ in
\eqref{eq:termination} and the total number of cost functional evaluations
50,000 whichever is reached first.

For generating $N$-collection of samples discussed in Section~\ref{sec:step_0}
we use $N = 10000$ and $N_{c,max} = 8$. The set of sample solutions
$\mC(10000)$ is pre-computed using a generator of uniformly distributed
random numbers. Therefore, each sample $\bar \sigma_i(x)$ ``contains''
from one to eight ``defective'' areas with $\sigma_d = 0.4$.
Each area is located randomly within domain $\Omega$ and represented
by a circle of randomly chosen radius $0 < r \leq 0.3 R$. Also, we
fix the number of samples $N_s$ to 10 for all numerical experiments
shown in this paper, see \cite{ArbicMS2020} for more details.
The entire preparation stage also involves computations for electrical
currents (measurements) using \eqref{eq:el_current}--\eqref{eq:forward} one
time per sample solution with total computing time equivalent roughly to
10,000 cost functional evaluations. For example, this phase is
completed in 10.3~h on a computer with a 26-core Intel Xeon Cascade
Lake Refresh Gold 6230R processor running at 2.1~GHz using 256~GB of RAM
under MS Windows Server 2019 by utilizing only one core. As all samples are
independent of each other, the associated computations are highly parallelizable;
this adds more to the computational efficiency of our new approach.

\subsection{Framework Validation}
\label{sec:model_valid}

\begin{figure*}[!htb]
  \begin{center}
  \mbox{
  \subfigure[model~\#1: $\sigma_{true}(x)$]
    {\includegraphics[width=0.33\textwidth]{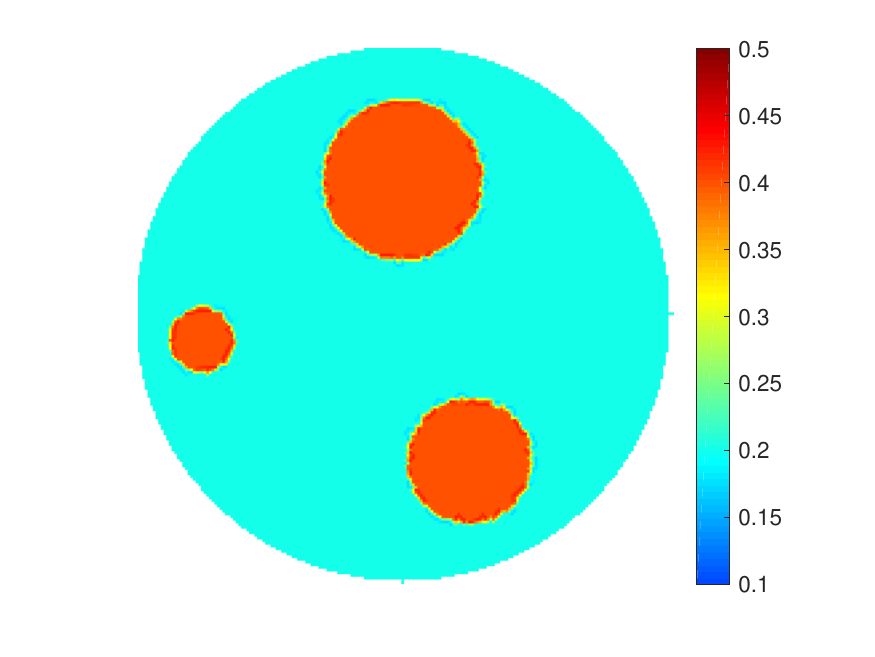}}
  \subfigure[best sample $\bar \sigma_1(x)$]
    {\includegraphics[width=0.33\textwidth]{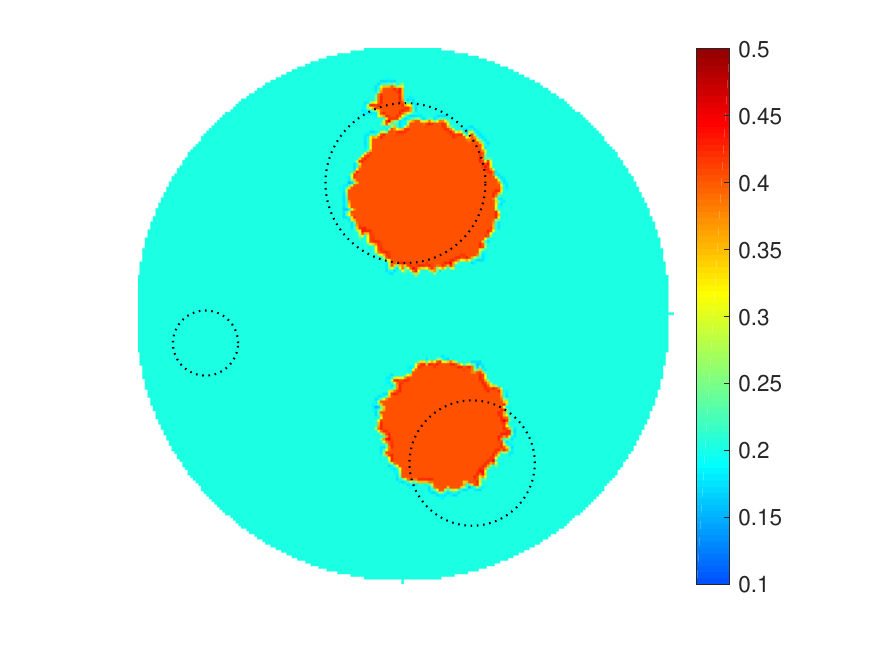}}
  \subfigure[solution $\sigma^0(x)$ after Step~1]
    {\includegraphics[width=0.33\textwidth]{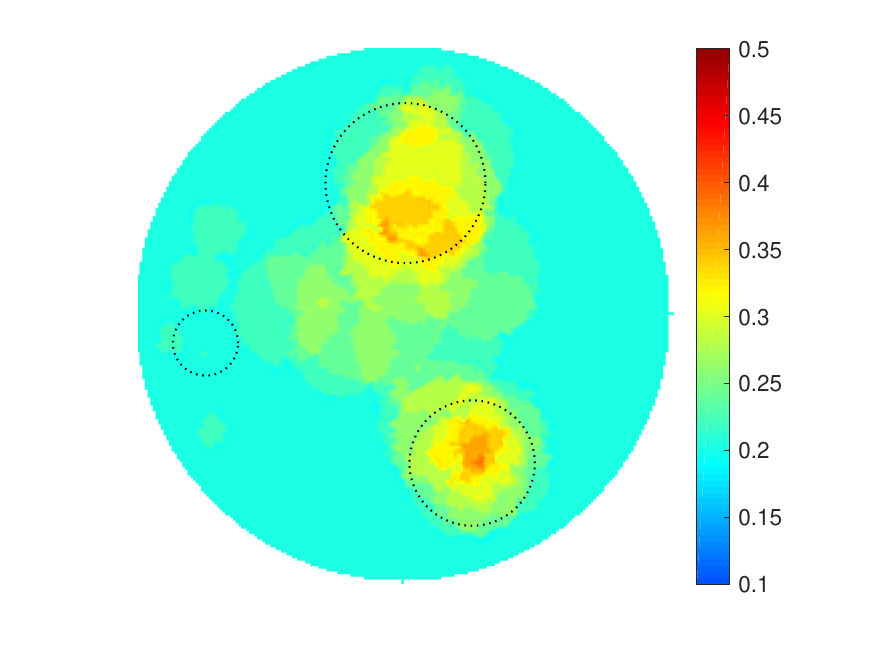}}}
  \end{center}
  \caption{Model~\#1. (a)~True electrical conductivity $\sigma_{true}(x)$.
    Solutions by (b)~best sample $\bar \sigma_1(x)$ from initial basis $\mB^0$,
    and (c)~complete initial basis $\mB^0$ approximation $\sigma^0(x)$ obtained
    after Step~1.}
  \label{fig:model_8a}
\end{figure*}

\begin{figure*}[!htb]
  \begin{center}
  \mbox{
  \subfigure[$\hat \sigma(x)$: no circles added]
    {\includegraphics[width=0.33\textwidth]{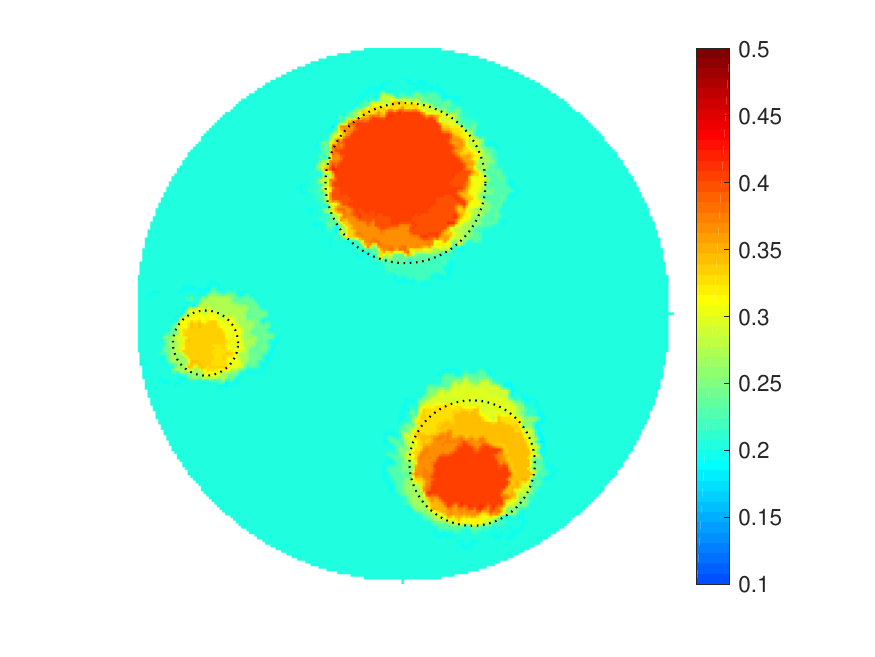}}
  \subfigure[$\hat \sigma(x)$: 8 circles/sample]
    {\includegraphics[width=0.33\textwidth]{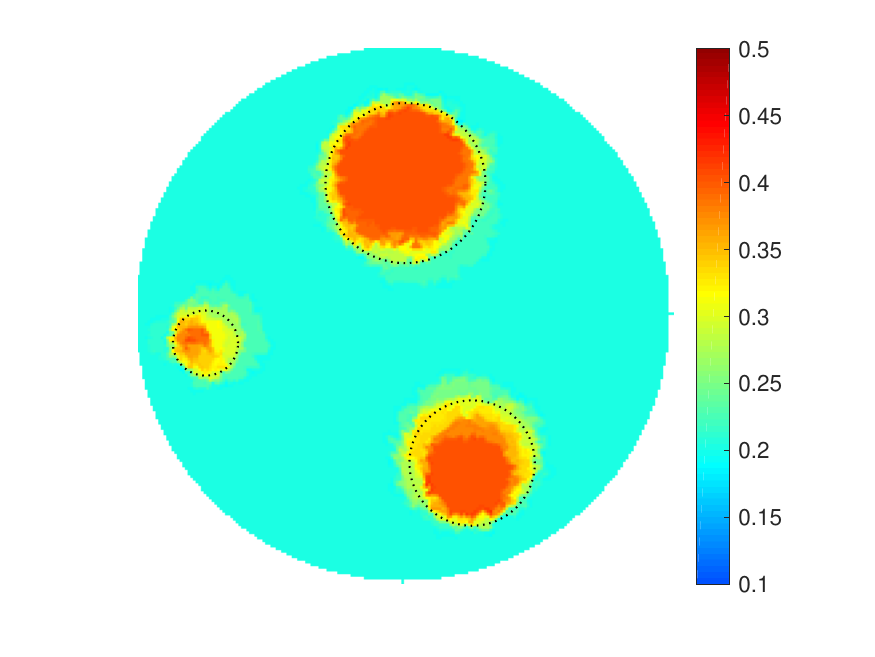}}
  \subfigure[effect of re-setting $N_c^i$]
    {\includegraphics[width=0.33\textwidth]{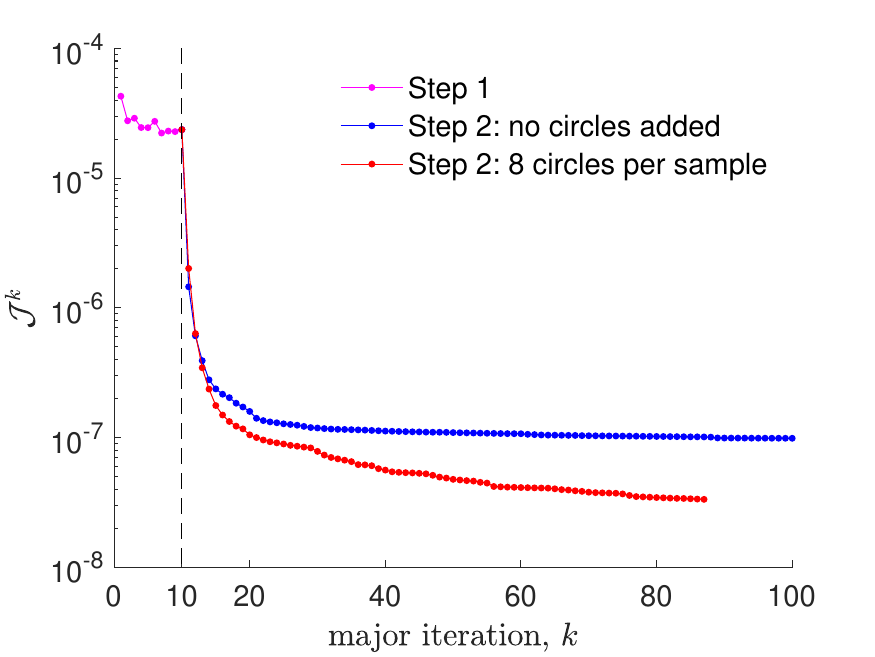}}}
  \mbox{
  \subfigure[$\alpha$: no circles added]
    {\includegraphics[width=0.33\textwidth]{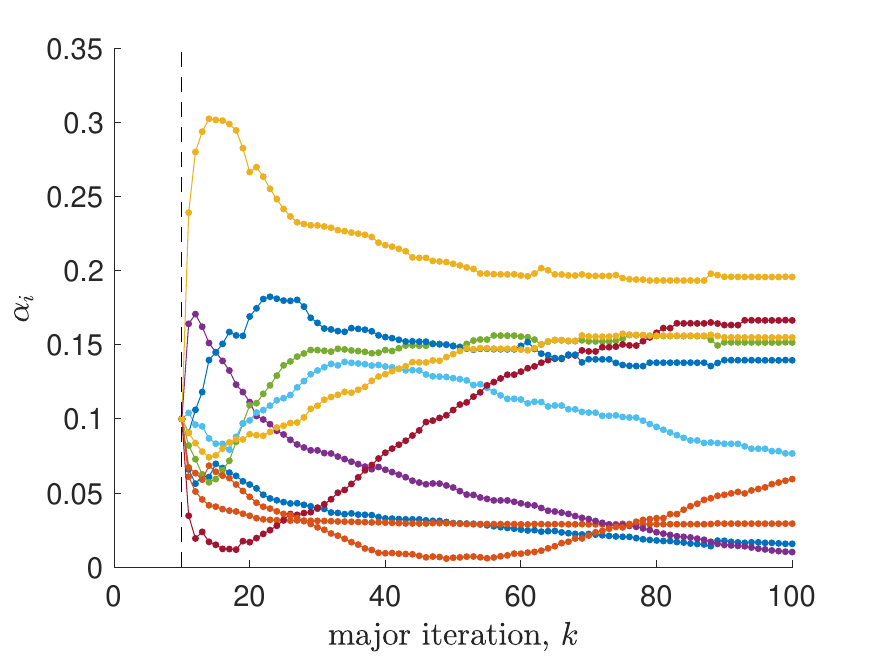}}
  \subfigure[$\alpha$: 8 circles/sample]
    {\includegraphics[width=0.33\textwidth]{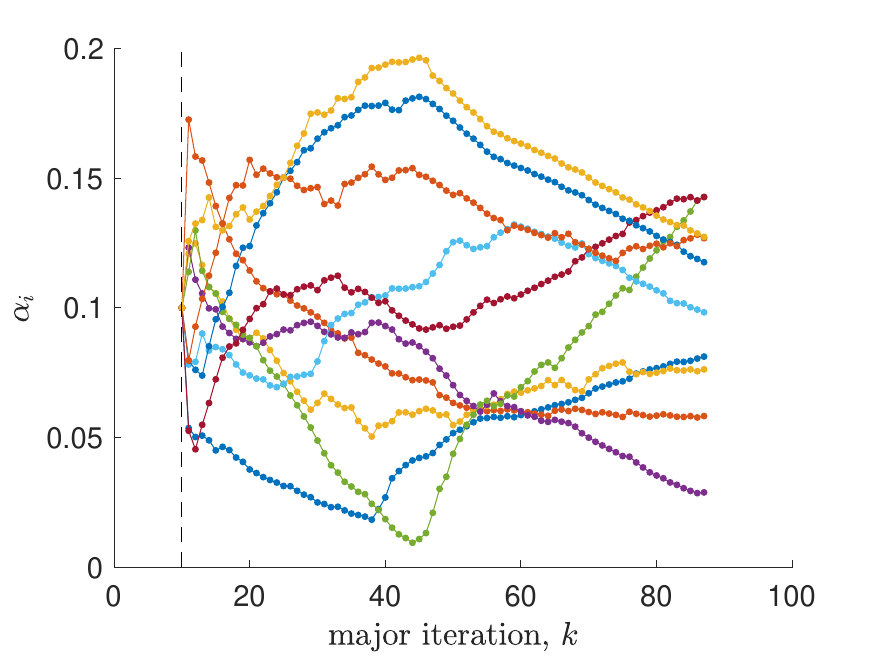}}
  \subfigure[solution error]{\includegraphics[width=0.33\textwidth]
    {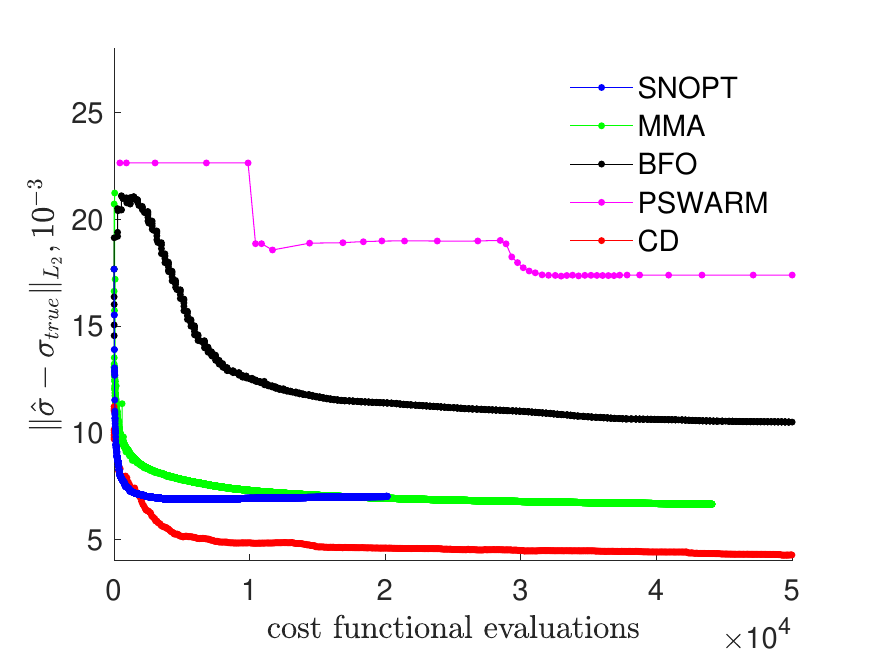}}}
  \end{center}
  \caption{Model~\#1. (a,b)~Solution images after Step~2 with (d,e)~history
    of control $\alpha$ changes for cases with (a,d)~no circles added
    and (b,e)~after adding new circles to each sample enforcing $N_c^i = 8$.
    (c)~Cost functionals $\mJ^k$ as a function of major iteration count
    $k$ evaluated at Step~1 (pink dots) and Step~2 with no circles added
    to samples (blue dots) and 8 circles per sample (red dots).
    (f)~Solution error $\| \hat \sigma - \sigma_{true} \|_{L_2}$
    as a function of number of cost functional evaluations for results
    obtained by gradient-based SNOPT (blue) and MMA (green),
    derivative-free BFO (black) and PSWARM (pink) all using PCA compared
    with CD (red) method with enforced $N_c^i = 8$ condition.}
  \label{fig:model_8b}
\end{figure*}

To begin checking the performance of the proposed optimization
framework, we created our (benchmark) model~\#1 to mimic a case when
a medium contains three areas of different sizes and circular
shapes to appear defective as seen in Figure~\ref{fig:model_8a}(a).
To prove the superior performance of the proposed algorithm
supplied with the customized CD method for optimization, we refer
first to Figure~\ref{fig:model_8a}(b,c). The center image (b) shows
the best, in terms of the measurement fit, sample solution found in
the collection $\mC(10000)$ for model~\#1 and placed into the initial
basis $\mB^0$ as $\bar \sigma_1(x)$. The right image (c) shows the
solution $\sigma^0(x)$ obtained after Step~1 optimization procedure.
Next, Figure~\ref{fig:model_8b}(a,b) depicts the solutions obtained
after Step~2 respectively with no circles and after adding new ones
to each sample. It is noticeable that the structure of the best sample
solution $\bar \sigma_1(x)$ is very far from $\sigma_{true}(x)$, and
$\sigma^0(x)$, derived from the initial basis $\mB^0$, poorly approximates
model~\#1. However, the proposed sample-based parameterization enables
the framework to accurately locate all defective regions, including the
smallest one. As the next step, we investigate the effect of re-setting parameters $N_c^i$
as discussed in Section~\ref{sec:step_1}. Figure~\ref{fig:model_8b}(c) shows
the progress during Step~1 (first 10 major iterations, in pink) and Step~2
with no circles added (in blue) and after adding new circles to each
sample (in red) so that $N_c^i = 8, \ i = 1, \ldots, 10$. An interesting
observation could be made from Figure~\ref{fig:model_8b}(d,e) presenting
the history of control $\alpha$ changes during Step~2 optimization. The
comparison reveals more sensitivity (wider oscillation spans for all
$\alpha$'s) in the second case related to the increased size of control
space $\mP$. Due to obviously better performance for that case, we
use the same strategy, by setting $N_c^i = N_{c,max}$, in all numerical
experiments presented in this paper unless otherwise stated.


We also use model~\#1 to compare the performance of the proposed
framework with those observed after applying various methods, namely
\begin{itemize}
  \item gradient-based sequential quadratic and sequential convex programming
    algorithms by, respectively, Sparse Nonlinear OPTimizer ({\tt SNOPT})
    \cite{SNOPTManual} and method of moving asymptotes ({\tt MMA})
    \cite{Svanberg1987},
  \item derivative-free brute force optimization~(BFO) method from the {\tt BFO~2.0} package
    \cite{BFO2017,BFO2020} and the dual pattern and particle swarm algorithm
    Particle SWARM ({\tt PSWARM}) \cite{PSWARM2007} implemented through the {\tt OPTI} Toolbox
    for nonlinear optimization problems \cite{OPTI2012} and
  \item our CD method customized to a predefined order of controls as described
    in Section~\ref{sec:solution}.
\end{itemize}
To compare the quality of the obtained solutions both gradient-based and
derivative-free methods (excluding our CD algorithm) are supplied with
principal component analysis (PCA) techniques for control space reduction as described in detail in
\cite{KoolmanBukshtynov2021,Bukshtynov2015,VolkovBukshtynov2018}. The PCA is
performed on the same collection of samples $\mC(10000)$ and with the same
number of principal components, 250 (preserving about 90\% of the ``energy''
in the full set of basis vectors), as the total number of controls used in
the CD method.
\begin{figure*}[!htb]
  \begin{center}
  \mbox{
  \subfigure[$\hat \sigma(x)$ by SNOPT]
    {\includegraphics[width=0.33\textwidth]{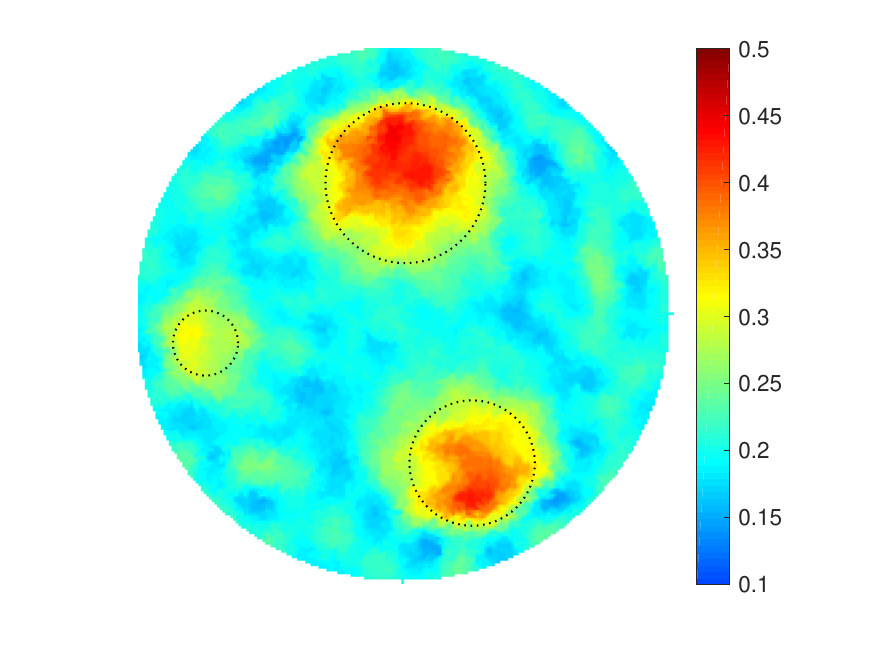}}
  \subfigure[$\hat \sigma(x)$ by MMA]
    {\includegraphics[width=0.33\textwidth]{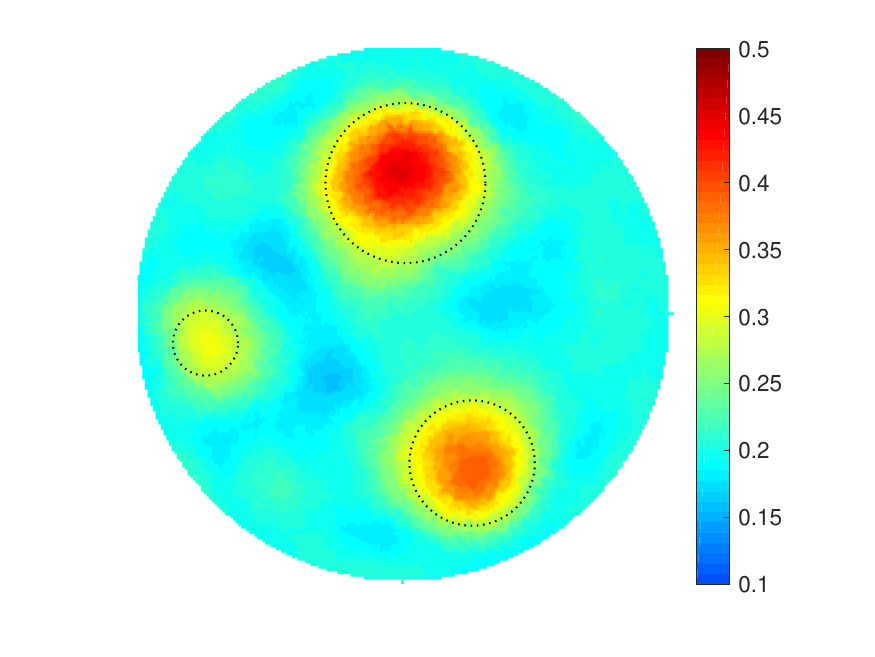}}
  \subfigure[$\hat \sigma(x)$ distribution]
    {\includegraphics[width=0.33\textwidth]{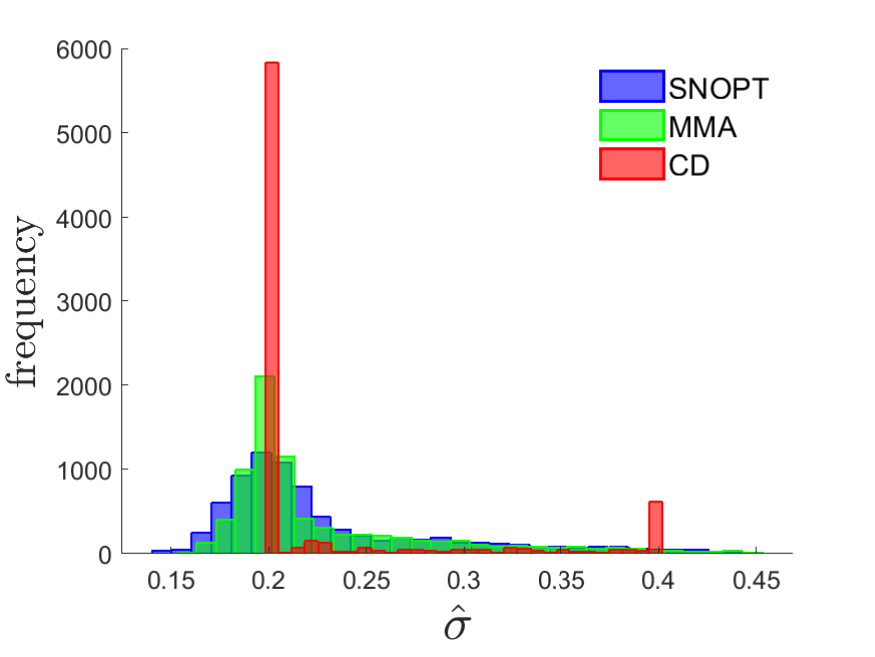}}}
  \mbox{
  \subfigure[$\hat \sigma(x)$ by BFO]
    {\includegraphics[width=0.33\textwidth]{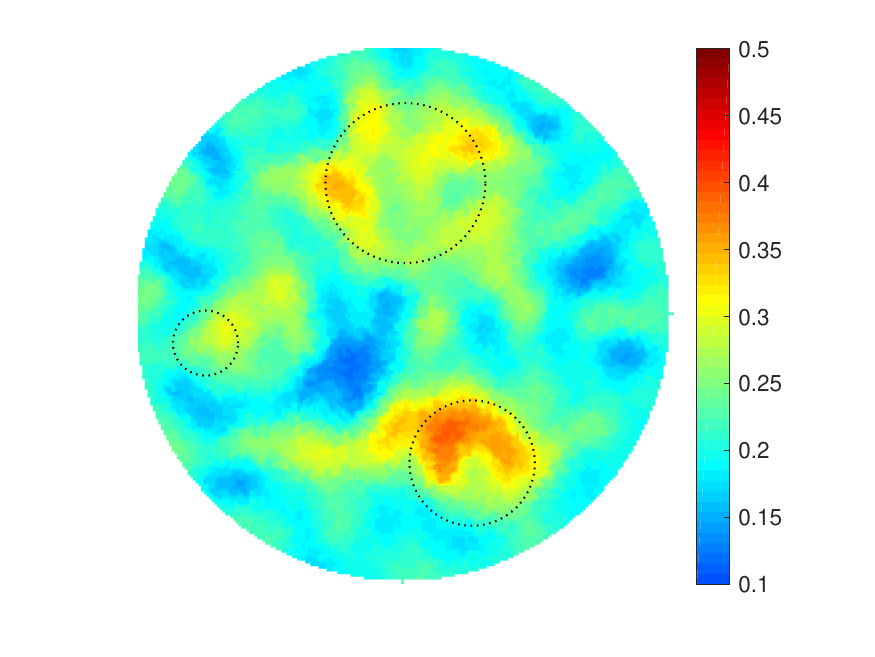}}
  \subfigure[$\hat \sigma(x)$ by PSWARM]
    {\includegraphics[width=0.33\textwidth]{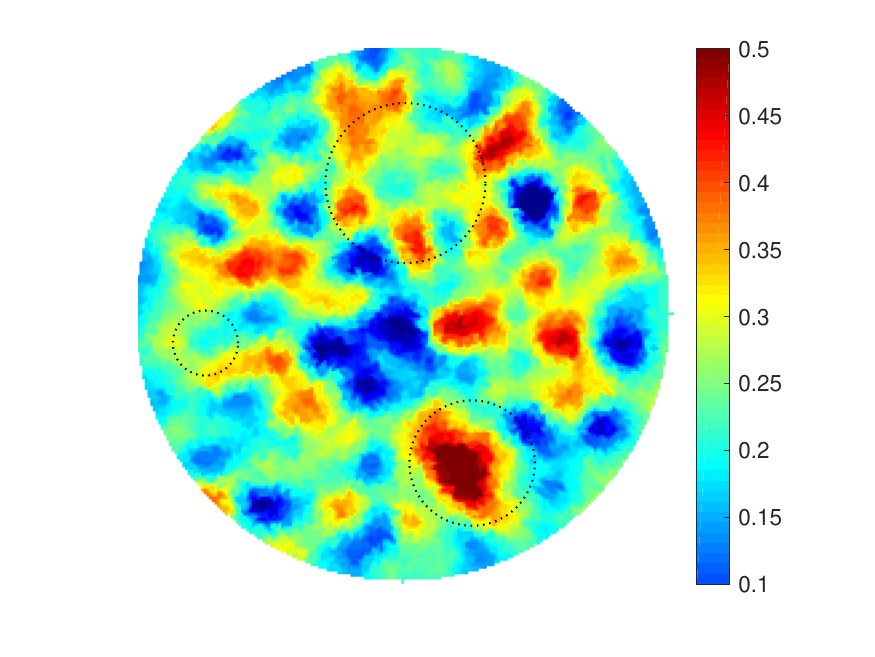}}
  \subfigure[$\hat \sigma(x)$ distribution]
    {\includegraphics[width=0.33\textwidth]{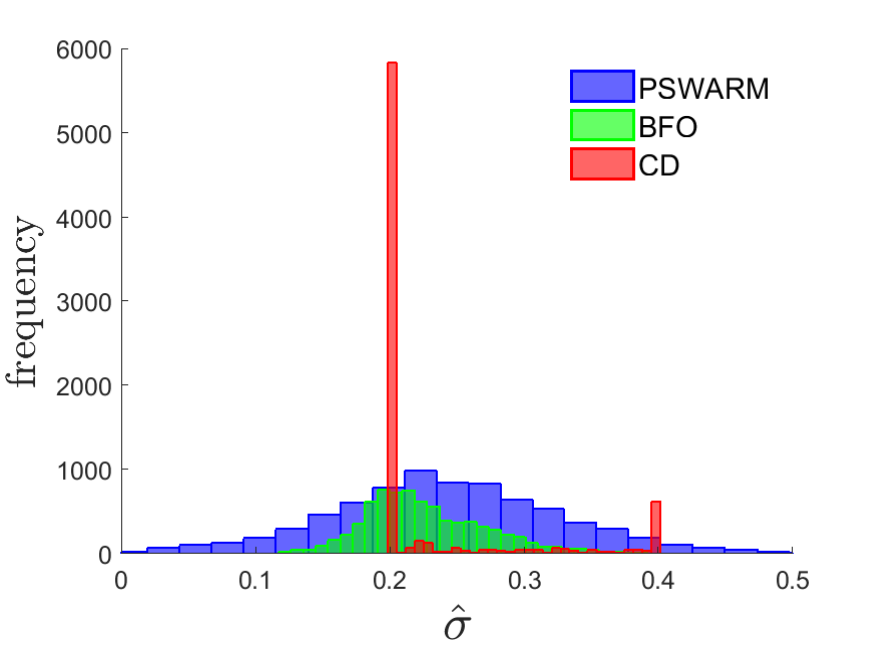}}}
  \end{center}
  \caption{Model~\#1. (a,b,d,e)~Solutions images obtained by (a)~SNOPT, (b)~MMA, (d)~BFO,
    and (e)~PSWARM. (c,f)~Histograms for solutions obtained by (c)~SNOPT and MMA,
    (f)~BFO and PSWARM compared with the CD method.}
  \label{fig:model_8c}
\end{figure*}

Figure~\ref{fig:model_8c}(a,b,d,e) shows the images obtained, respectively,
by SNOPT, MMA, BFO, and PSWARM approaches with combined accuracy results
provided in Figure~\ref{fig:model_8b}(f). To make the superior performance of
the proposed customized CD method even more distinctive, all four approaches
used termination tolerance $\epsilon = 10^{-9}$ in contrast to $10^{-4}$ in
the CD's case. As seen in Figure~\ref{fig:model_8b}(f), in general, gradient-based
approaches are able to provide images of better qualities. Both derivative-free
approaches, as expected, showed much slower convergence terminated when the
total number of cost functional evaluations reached 50,000. Although our derivative-free
CD approach is also terminated with the same condition, it arrives at a solution
with better quality than is seen in the gradient-based methods. To add more,
only our new CD method results in quality close to the desired binary distribution
(refer to Figure~\ref{fig:model_8c}(c,f) for analysis by histograms).

\begin{figure*}[!htb]
  \begin{center}
  \mbox{
  \subfigure[$\sigma_{true}(x), \sigma_d = 0.3$]
    {\includegraphics[width=0.25\textwidth]{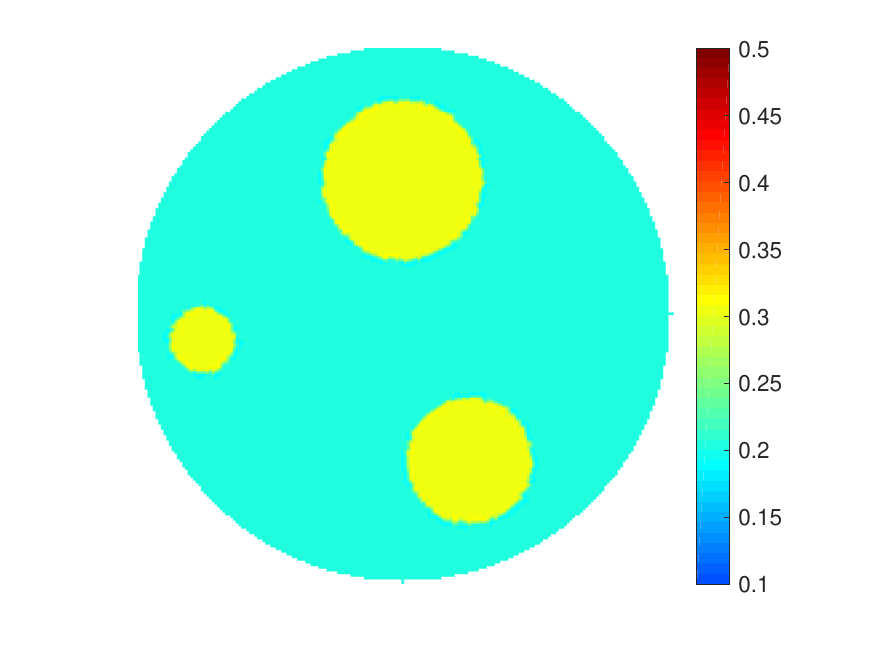}}
  \subfigure[$\hat \sigma(x), \sigma_d = 0.3$]
    {\includegraphics[width=0.25\textwidth]{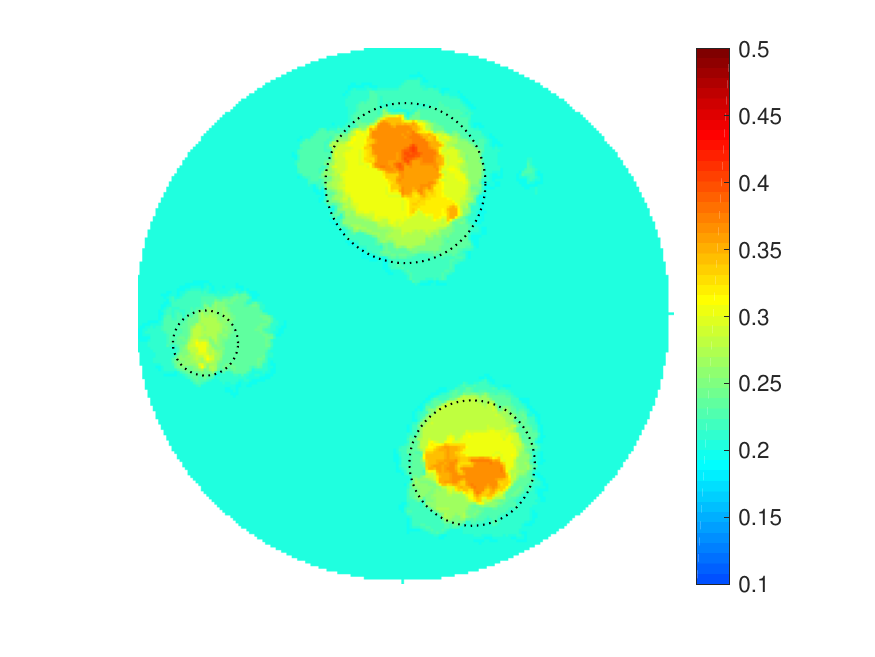}}
  \subfigure[$\sigma_{true}(x), \sigma_d = 0.5$]
    {\includegraphics[width=0.25\textwidth]{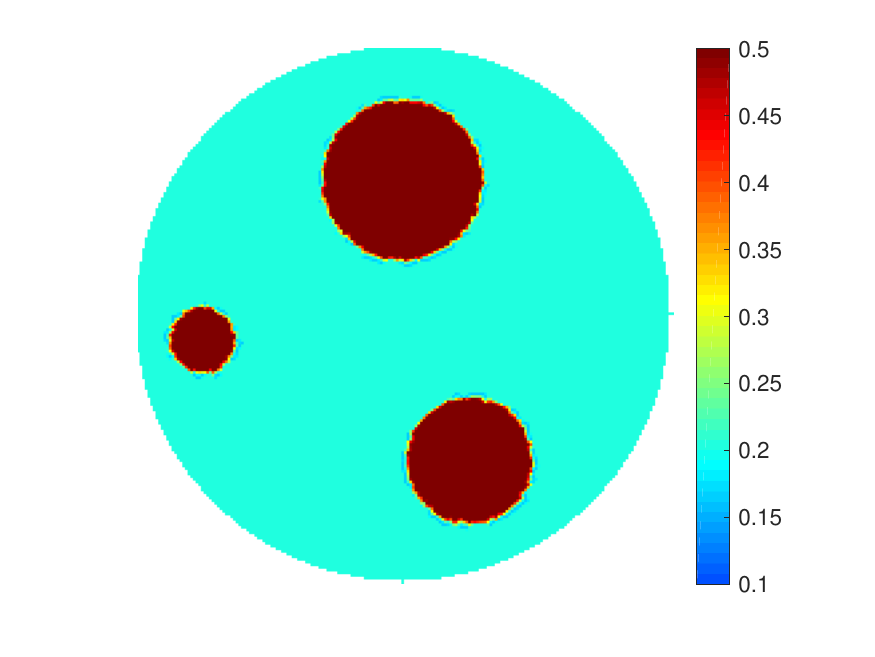}}
  \subfigure[$\hat \sigma(x), \sigma_d = 0.5$]
    {\includegraphics[width=0.25\textwidth]{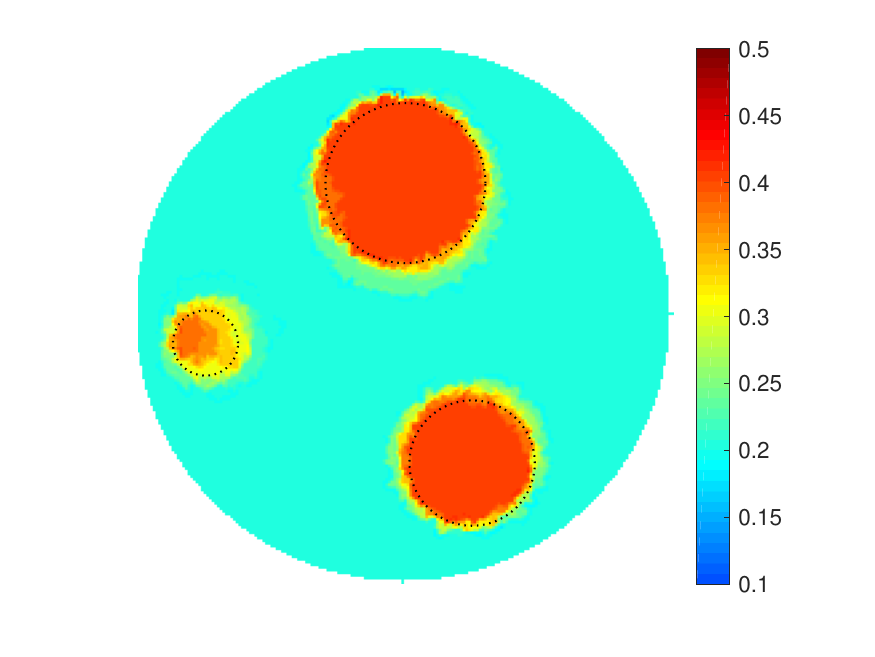}}}
  \end{center}
  \caption{Model~\#1. (a,c)~True electrical conductivity $\sigma_{true}(x)$
    and (b,d)~respective solution images obtained by the derivative-free
    customized CD method when (a,b)~$\sigma_d = 0.3$ and (c,d)~$\sigma_d = 0.5$.}
  \label{fig:model_8d}
\end{figure*}

Finally, we would like to address the issue of {\it a priori} knowledge of $\sigma_d$
and $\sigma_n$ used in \eqref{eq:smpl_par} as approximations for conductivities,
respectively, of defective region $\Omega_d$ and the rest of the domain $\Omega_n$.
The solution formulation given by \eqref{eq:sigma_main} suggests that values $\sigma_d$
and $\sigma_n$ in \eqref{eq:sigma_true} may vary from those used in the sampling procedure
by \eqref{eq:smpl_par}. We check this ability for reconstructing different conductivity
values by running two numerical experiments where $\sigma_d$ is set to 0.3 and 0.5
(instead of 0.4 used before). Figure~\ref{fig:model_8d} presents the results from both
cases proving the method's ability to detect the inclusions accurately. Further comparison
of images in Figures~\ref{fig:model_8d}(b,d) and \ref{fig:model_8b}(b) (varied vs.~predefined
conductivity) reveals the difference in the $\sigma$-values within the defective regions
and suggests the suitability of the proposed methodology for quantitative imaging. The entire
computational framework may be further modified to allow both $\sigma_d$ and $\sigma_n$ to be
randomly chosen at the expense of enlarging the control space. It may require adding gradients
combined with the derivative-free approach considered in this paper to keep the method
computationally efficient.

\subsection{Effect of Noise in Data}
\label{sec:noise}

Now we would like to address a well-known issue of the noise present
in the measurements due to improper electrode-medium contacts,
possible electrode misplacement, wire interference, etc. The effect
of noise has already been investigated by many researchers both
theoretically and within practical applications with suggested
approaches to mitigate its negative impact on the quality of images.
In this section, we compare the effect of noise in reconstructions
obtained by the gradient-based SNOPT and MMA (both with PCA) and
our proposed derivative-free customized CD method.
\begin{figure*}[!htb]
  \begin{center}
  \mbox{
  \subfigure[CD: 0.5\% noise]{\includegraphics[width=0.25\textwidth]
    {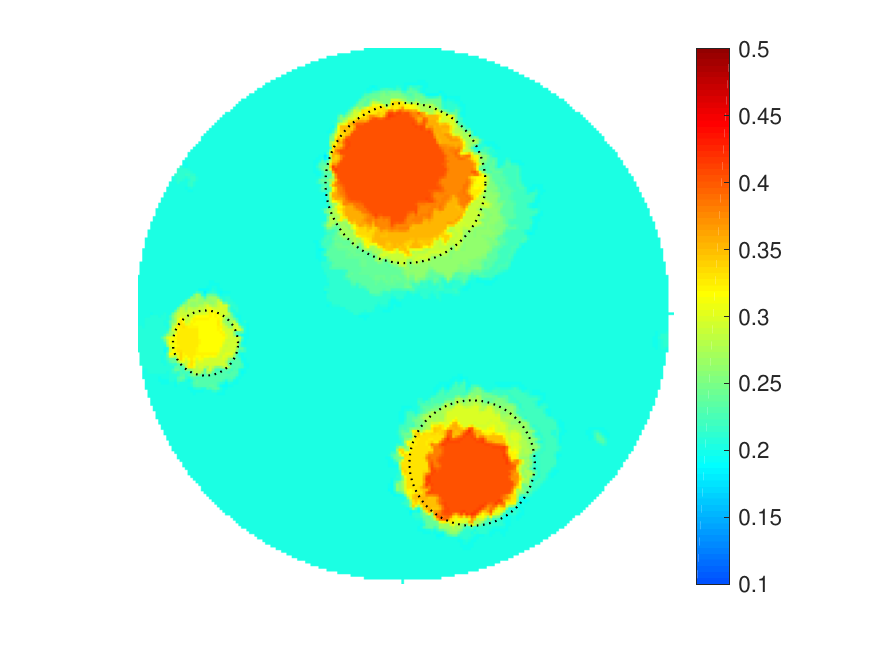}}
  \subfigure[1\% noise]{\includegraphics[width=0.25\textwidth]
    {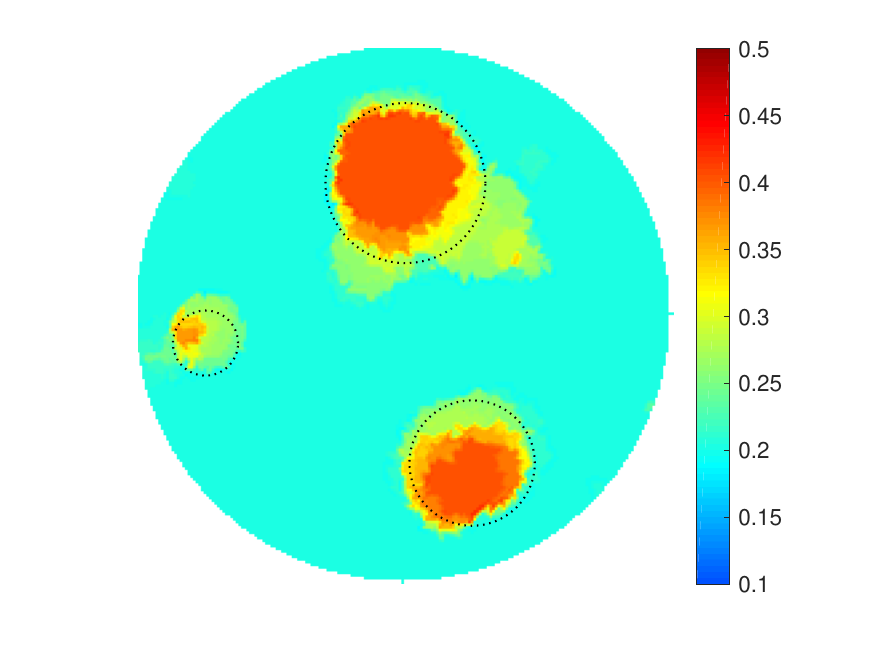}}
  \subfigure[2\% noise]{\includegraphics[width=0.25\textwidth]
    {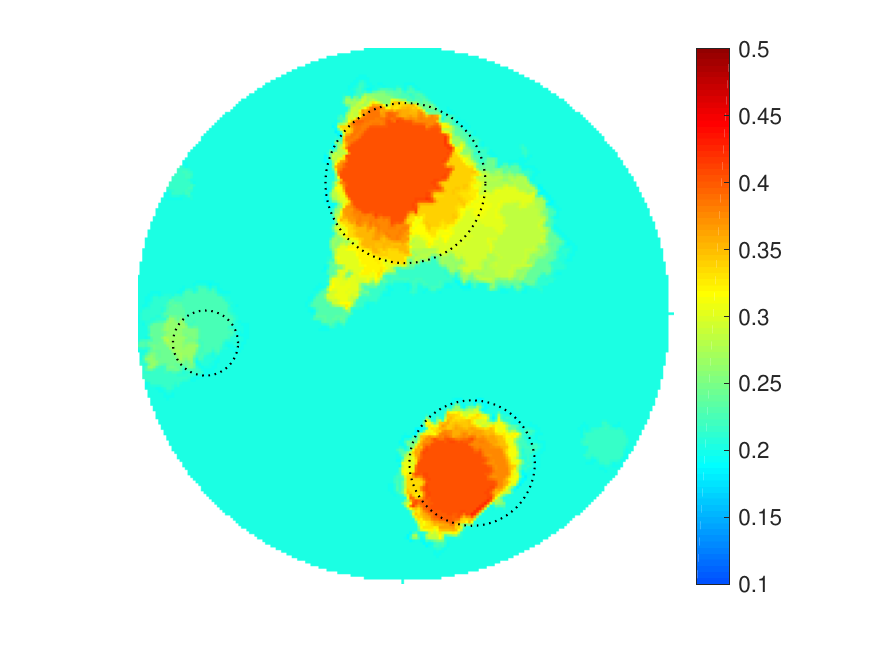}}
  \subfigure[5\% noise]{\includegraphics[width=0.25\textwidth]
    {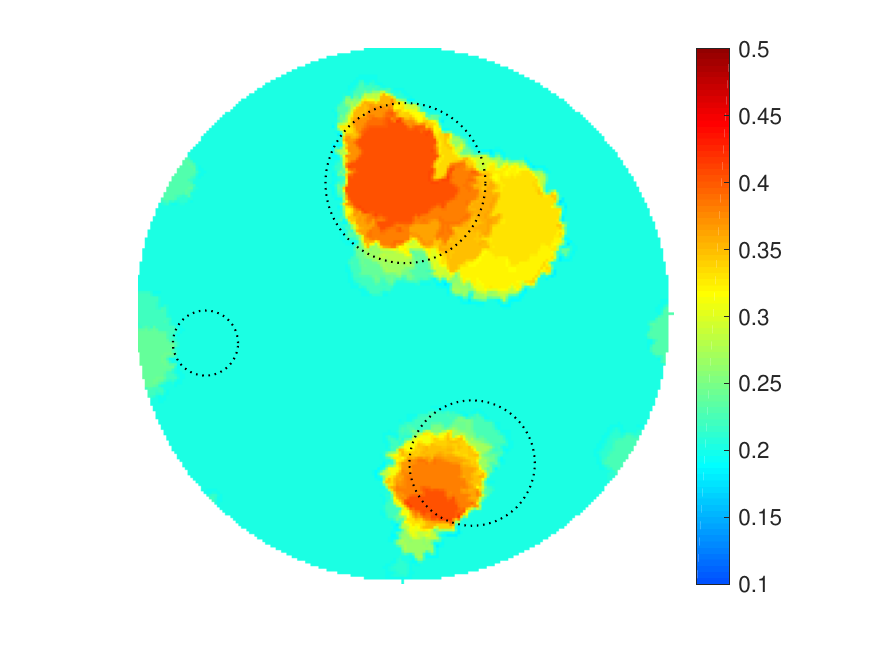}}}
  \mbox{
  \subfigure[SNOPT: 0.5\% noise]{\includegraphics[width=0.25\textwidth]
    {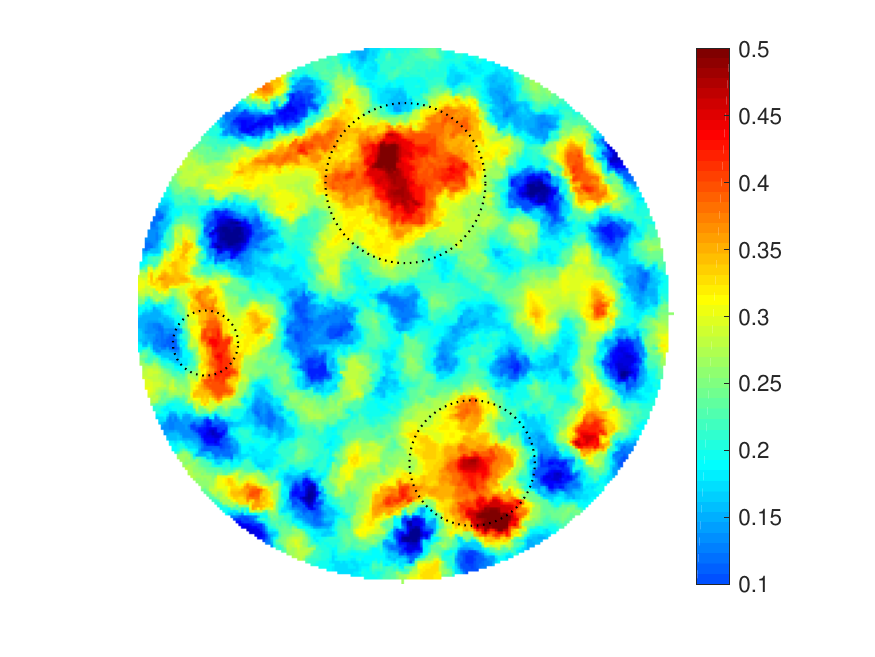}}
  \subfigure[1\% noise]{\includegraphics[width=0.25\textwidth]
    {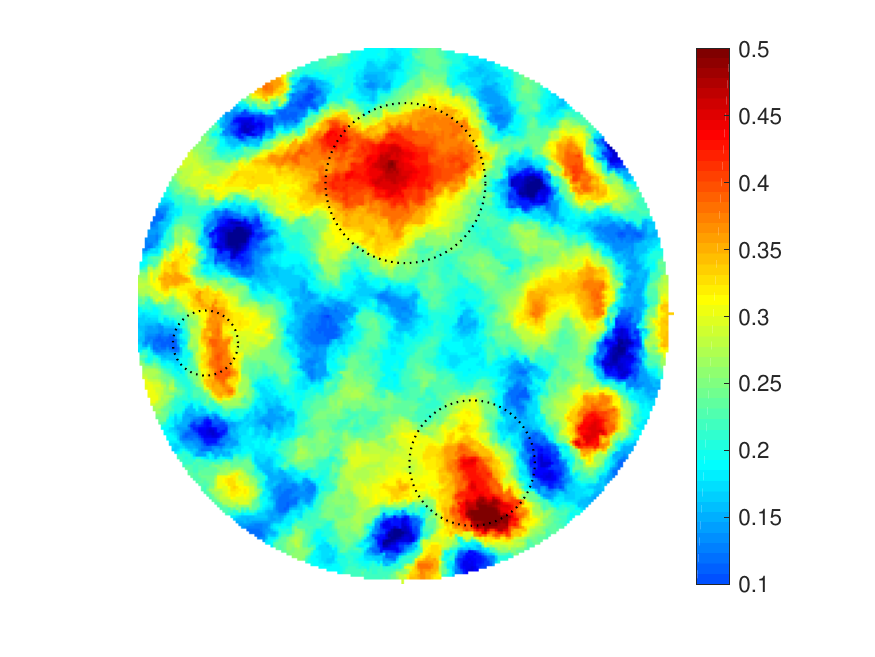}}
  \subfigure[2\% noise]{\includegraphics[width=0.25\textwidth]
    {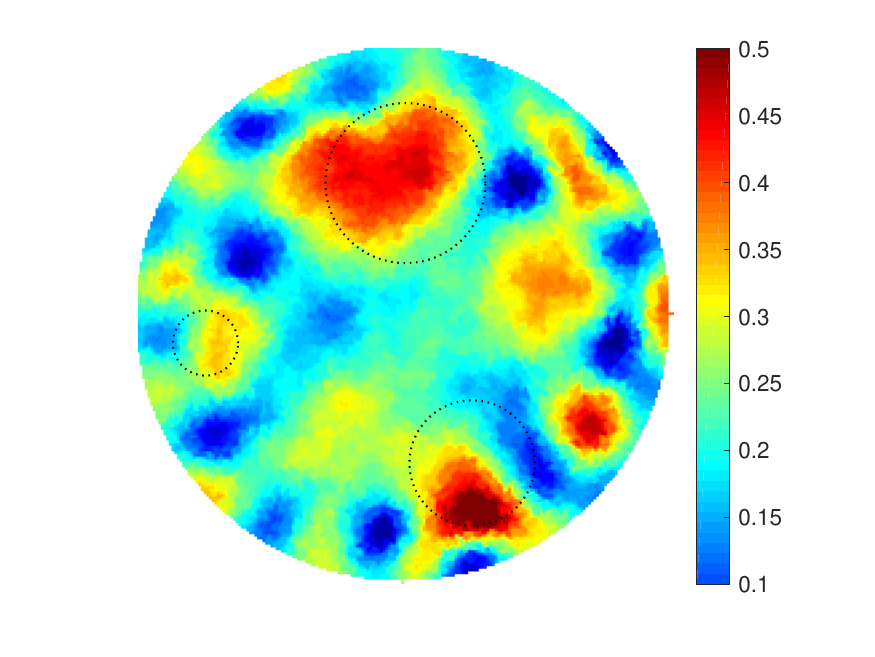}}
  \subfigure[5\% noise]{\includegraphics[width=0.25\textwidth]
    {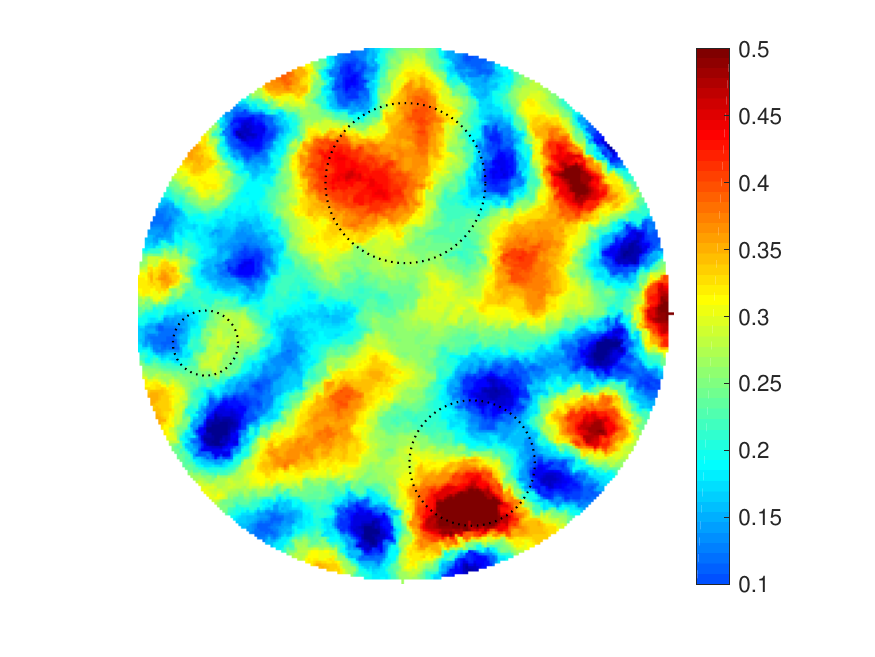}}}
  \mbox{
  \subfigure[MMA: 0.5\% noise]{\includegraphics[width=0.25\textwidth]
    {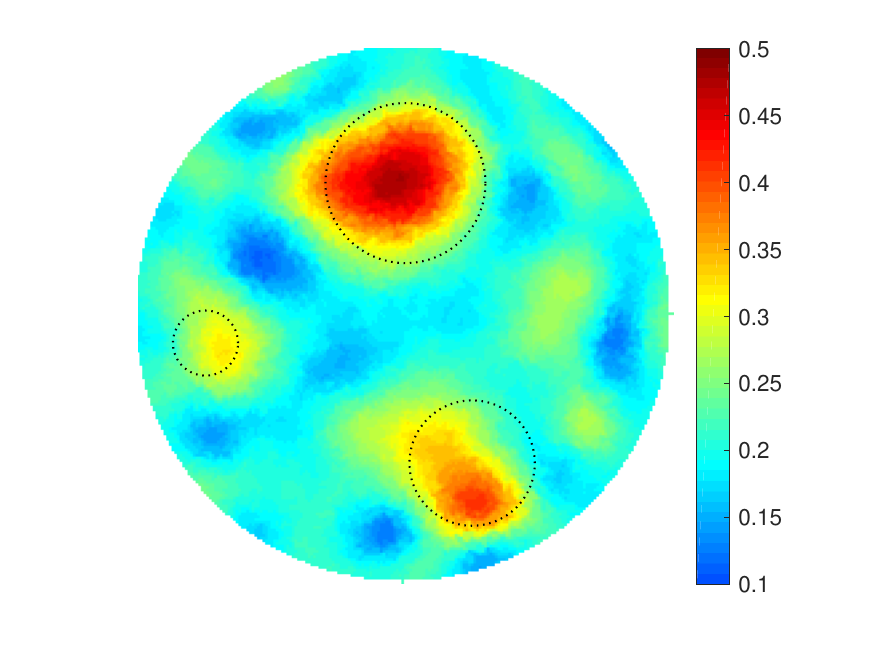}}
  \subfigure[1\% noise]{\includegraphics[width=0.25\textwidth]
    {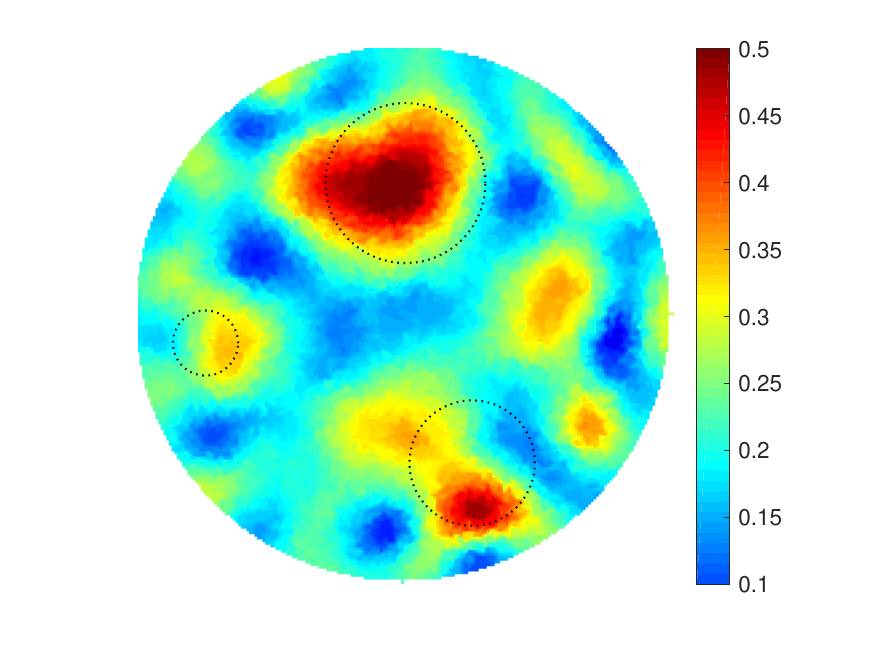}}
  \subfigure[2\% noise]{\includegraphics[width=0.25\textwidth]
    {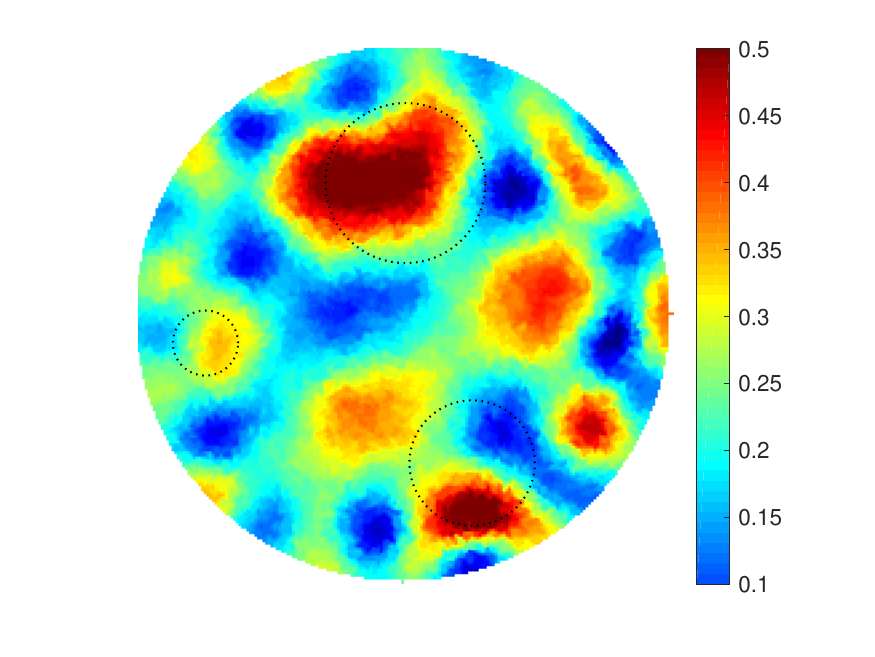}}
  \subfigure[5\% noise]{\includegraphics[width=0.25\textwidth]
    {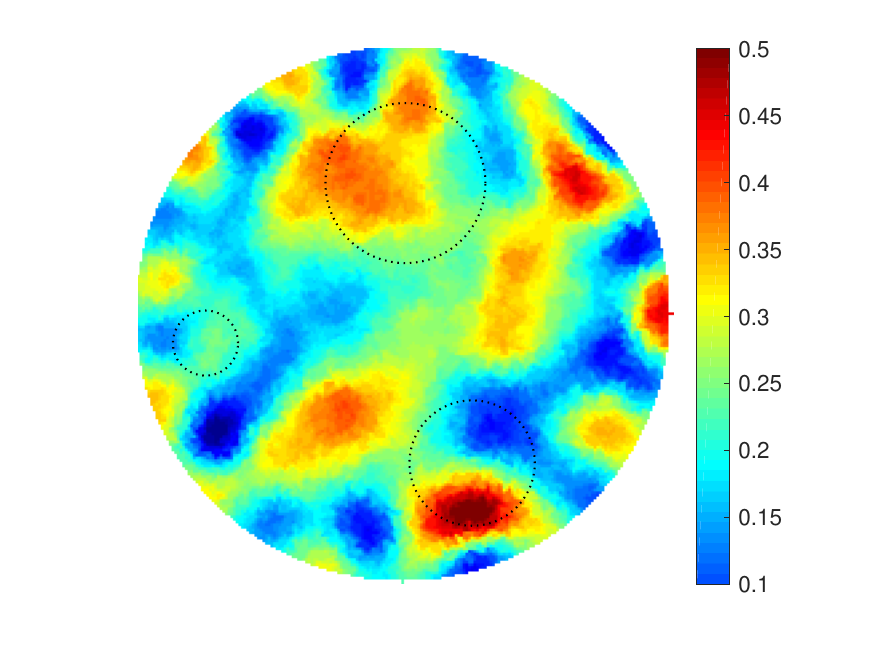}}}
  \end{center}
  \caption{Model~\#1. Solution images obtained by (a-d)~the
    proposed derivative-free customized CD, gradient-based
    (e-h)~SNOPT, and (i-l)~MMA methods when measurements are
    contaminated with (a,e,i)~0.5\%, (b,f,j)~1\%, (c,g,k)~2\%, and
    (d,h,l)~5\% noise.}
  \label{fig:model_8_noise_comp}
\end{figure*}

\begin{figure}[!htb]
  \begin{center}
  \includegraphics[width=0.75\textwidth]{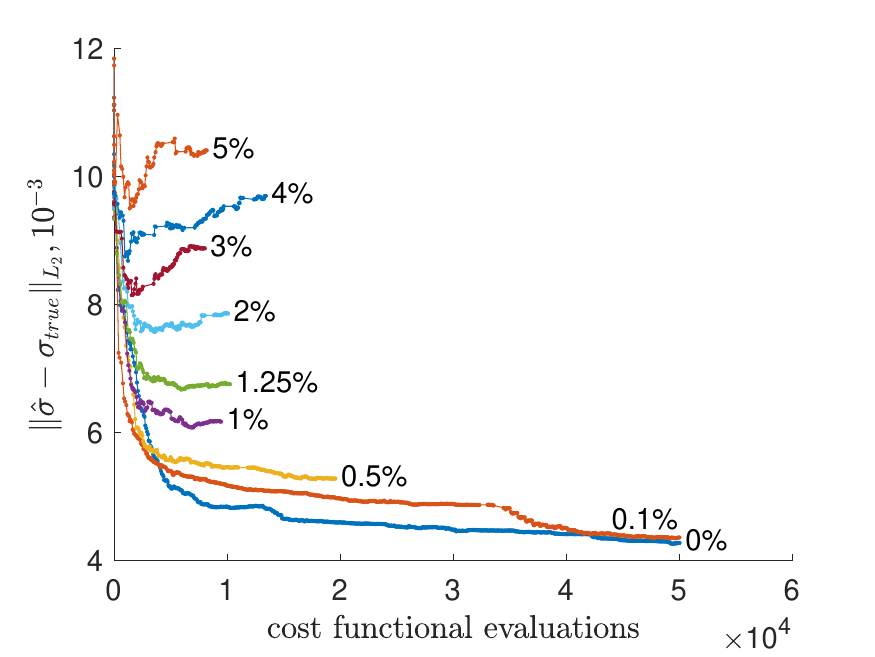}
  \end{center}
  \caption{Model~\#1. Solution error
    $\| \hat \sigma - \sigma_{true} \|_{L_2}$ as a function of
    the number of cost functional evaluations for results obtained
    by the proposed framework for various levels of noise.}
  \label{fig:model_8_noise_CD}
\end{figure}

In Figure~\ref{fig:model_8_noise_comp}, we revisit model~\#1 presented
first in Section~\ref{sec:model_valid} now with measurements
contaminated with 0.5\%, 1\%, 2\% and 5\% normally distributed noise.
As expected, we see that various levels of noise lead to oscillatory
instabilities in the images reconstructed by the gradient-based
approaches utilizing parameterization via PCA. This will obviously
result in multiple cases of false positive outcomes in screening procedures.
On the other hand, our new approach with sample-based parameterization shows
its stable ability to provide clear and accurate images with the
appearance of false negative results for small regions only with
noise higher than 2\%. Figure~\ref{fig:model_8_noise_CD} provides
the complete comparison of the solution error for results obtained
by our framework for various levels of noise between 0\% and 5\%.
We close this section by concluding on using 0.5\% noise for the
rest of the numerical experiments shown in this paper.

\subsection{Validation with Complicated Models}
\label{sec:models_complex}

In this section, we present results obtained using our new optimization
framework with sample-based parameterization applied to models with a
significantly increased level of complexity. The new algorithm has
already confirmed its ability to reconstruct accurately defective
regions of various size and at multiple locations. Therefore, the added
complications we are focusing on here are of small sizes and non-circular
shapes for those regions.

\begin{figure*}[!htb]
  \begin{center}
  \mbox{
  \subfigure[model~\#2: $\sigma_{true}(x)$]
    {\includegraphics[width=0.33\textwidth]{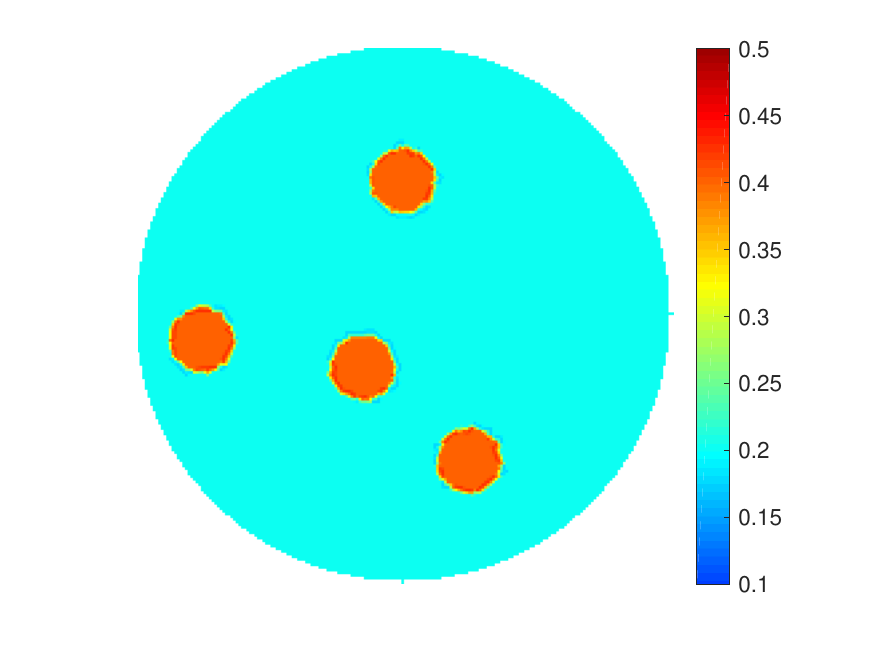}}
  \subfigure[CD: no noise]
    {\includegraphics[width=0.33\textwidth]{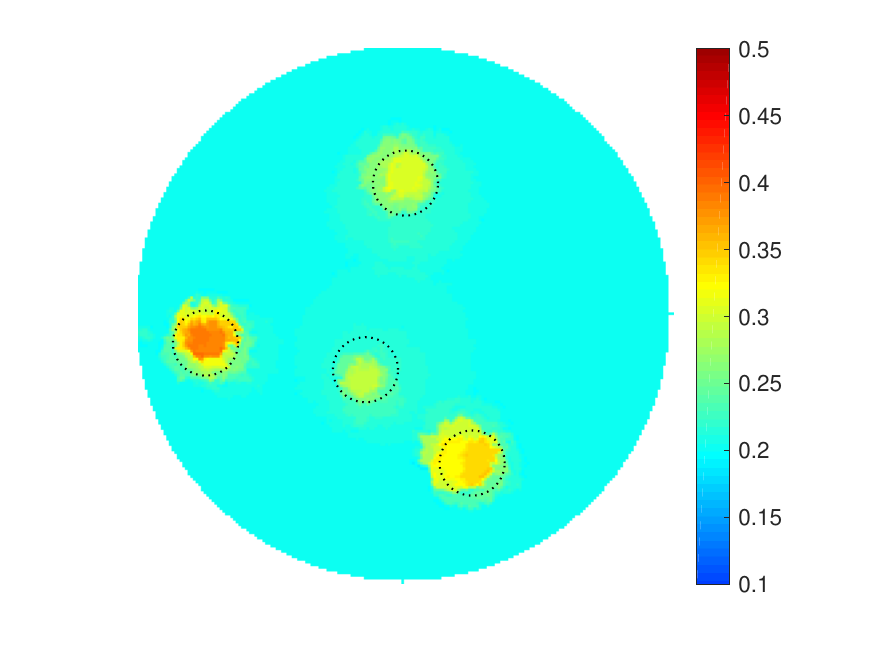}}
  \subfigure[CD: 0.5\% noise]
    {\includegraphics[width=0.33\textwidth]{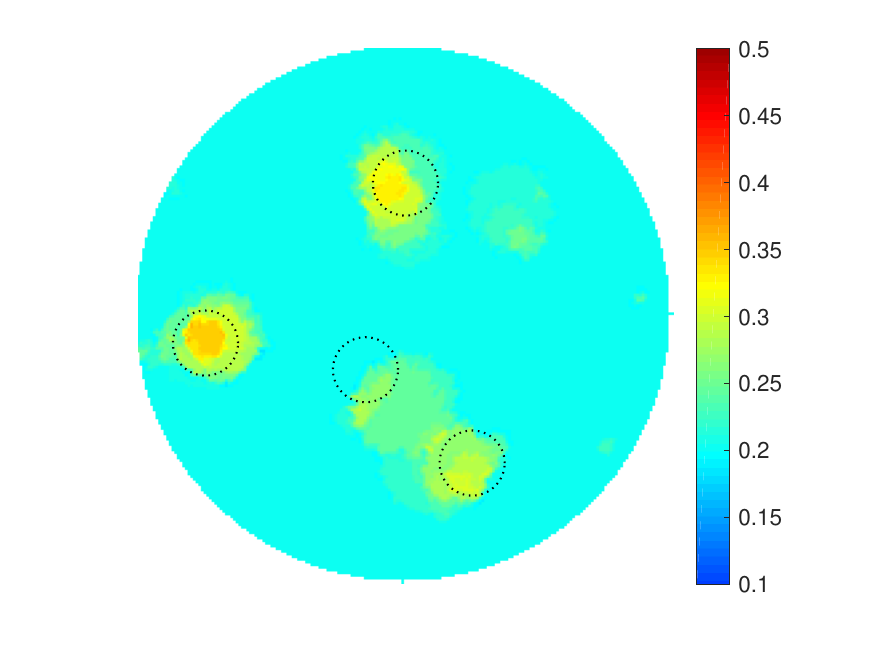}}}
  \mbox{
  \subfigure[solution error]{\includegraphics[width=0.33\textwidth]
    {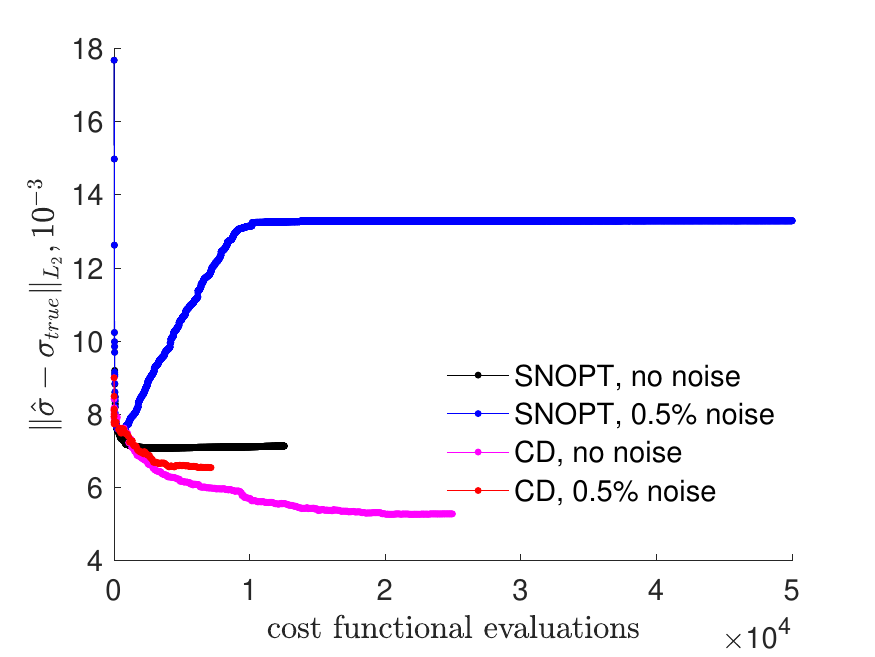}}
  \subfigure[SNOPT: no noise]{\includegraphics[width=0.33\textwidth]
    {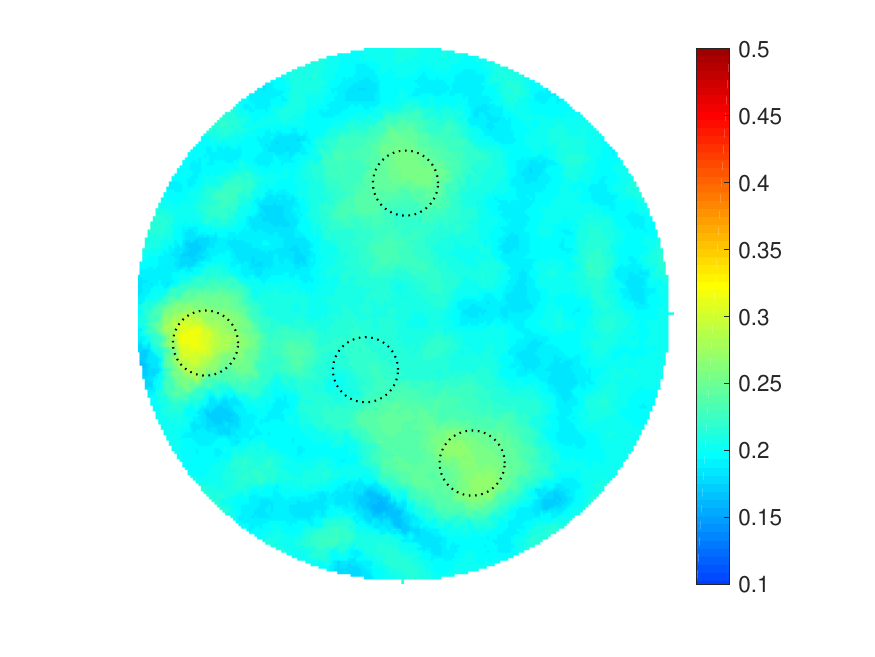}}
  \subfigure[SNOPT: 0.5\% noise]{\includegraphics[width=0.33\textwidth]
    {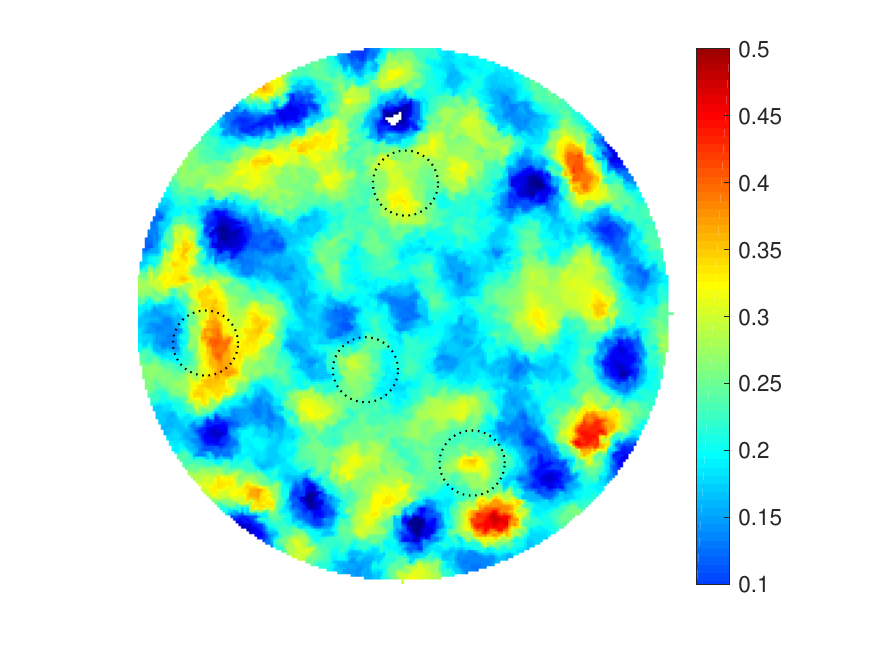}}}
  \end{center}
  \caption{Model~\#2. (a)~True electrical conductivity $\sigma_{true}(x)$.
    (b,c,e,f)~Solution images obtained by (b,c)~the proposed framework
    and (e,f)~the gradient-based SNOPT with PCA with (b,e)~no noise added
    and (c,f)~0.5\% noise in measurements. (d)~Solution error
    $\| \hat \sigma - \sigma_{true} \|_{L_2}$ as a function of the number
    of cost functional evaluations for images shown in (b,c,e,f).}
  \label{fig:model_3}
\end{figure*}

Our model \#2 is created to check the ability of the EIT techniques
coupled with our new approach to find the inclusions (defects) of very
small sizes. It mimics, for example, the application of EIT in medical practice
for recognizing cancer at early stages. The systematic analysis on the
distinguishability of two different conductivity distributions with
a specified precision is provided in \cite{Isaacson1986}. The electrical
conductivity $\sigma_{true}(x)$ is shown in Figure~\ref{fig:model_3}(a).
This model contains four circular-shaped defective regions all of the same size
as the smallest region in model \#1. The known complication comes from the
fact that the order of difference in measurements generated by this model
and ``fully normal medium'' ($\sigma(x) = \sigma_n, \ \forall x \in \Omega$)
is very close to the order of noise that appears naturally in the provided
data. In addition to this, small regions have a lower chance of being
detected if they are located closer to the center of the domain.

Figure~\ref{fig:model_3}(b-f) compares the results obtained by the
gradient-based SNOPT with PCA and our proposed derivative-free customized
CD methods without noise and with 0.5\% noise added to the measurements.
By analyzing images in Figure~\ref{fig:model_3}, we see that our approach
is able to provide more assistance in concluding on possible defects
in medium and help navigate for their search. When noise is negligible,
Figure~\ref{fig:model_3}(b), all four small spots are
distinguishable with accurately reconstructed shapes. Although adding
noise, as seen in Figure~\ref{fig:model_3}(c), brings more complexity to
image interpretation, it still keeps the possibility to help identify
defective regions. In fact we cannot claim the same for images
obtained by means of PCA and gradients, see Figure~\ref{fig:model_3}(e,f).
Figure~\ref{fig:model_3}(d) also proves computational efficiency of our
framework in comparison with the gradient-based method.

\begin{figure*}[!htb]
  \begin{center}
  \mbox{
  \subfigure[model~\#3: $\sigma_{true}(x)$]
    {\includegraphics[width=0.33\textwidth]{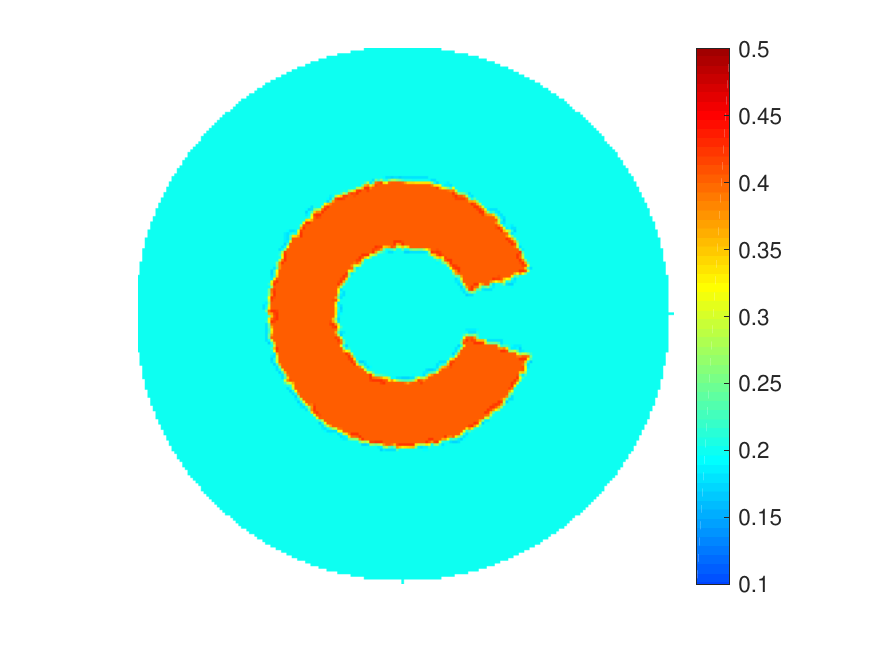}}
  \subfigure[CD: 8 circles, no noise]
    {\includegraphics[width=0.33\textwidth]{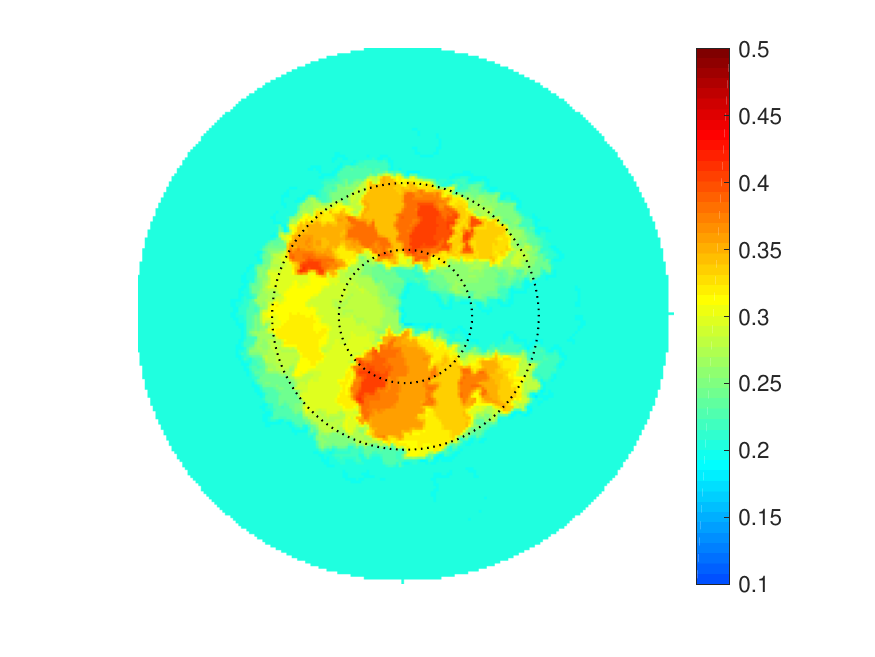}}
  \subfigure[CD: 8 circles, 0.5\% noise]
    {\includegraphics[width=0.33\textwidth]{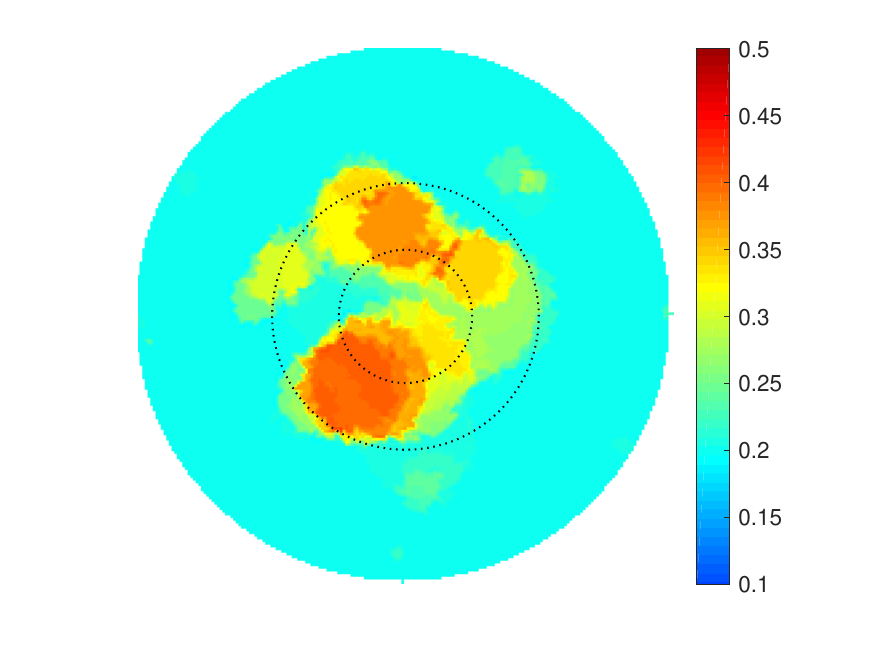}}}
  \mbox{
  \subfigure[SNOPT: no noise]{\includegraphics[width=0.33\textwidth]
    {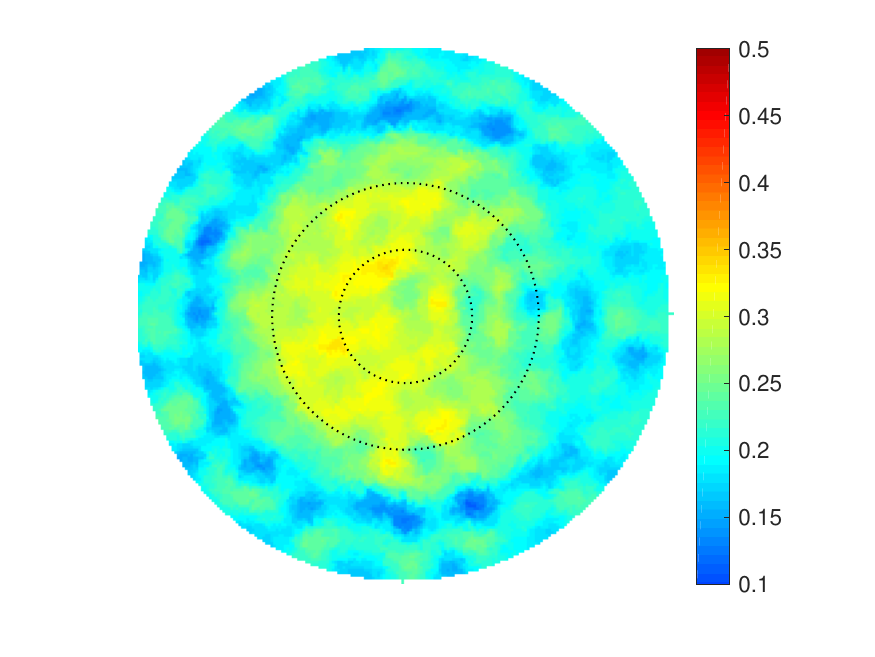}}
  \subfigure[CD: 20 circles, no noise]
    {\includegraphics[width=0.33\textwidth]{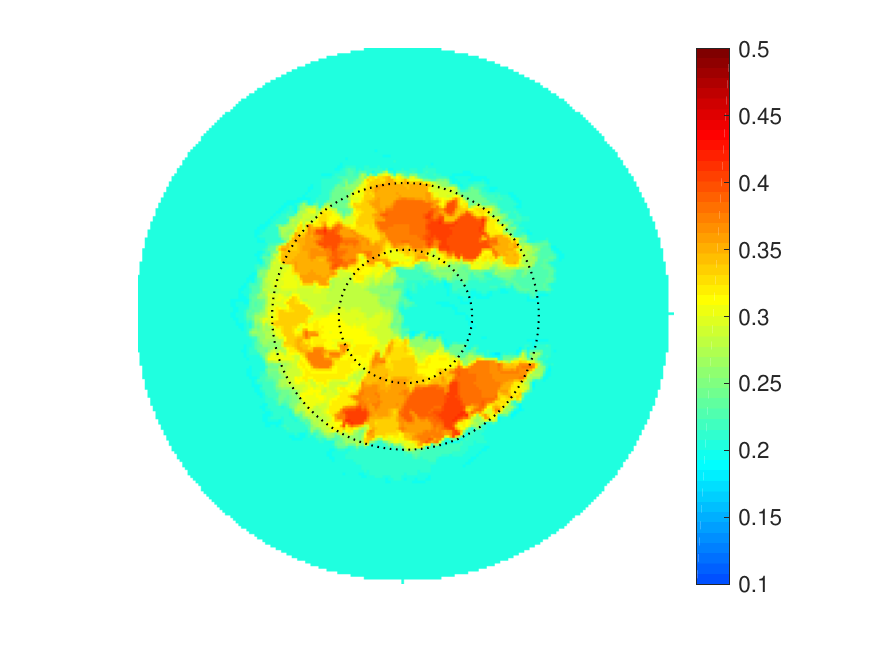}}
  \subfigure[CD: 20 circles, 0.5\% noise]
    {\includegraphics[width=0.33\textwidth]
    {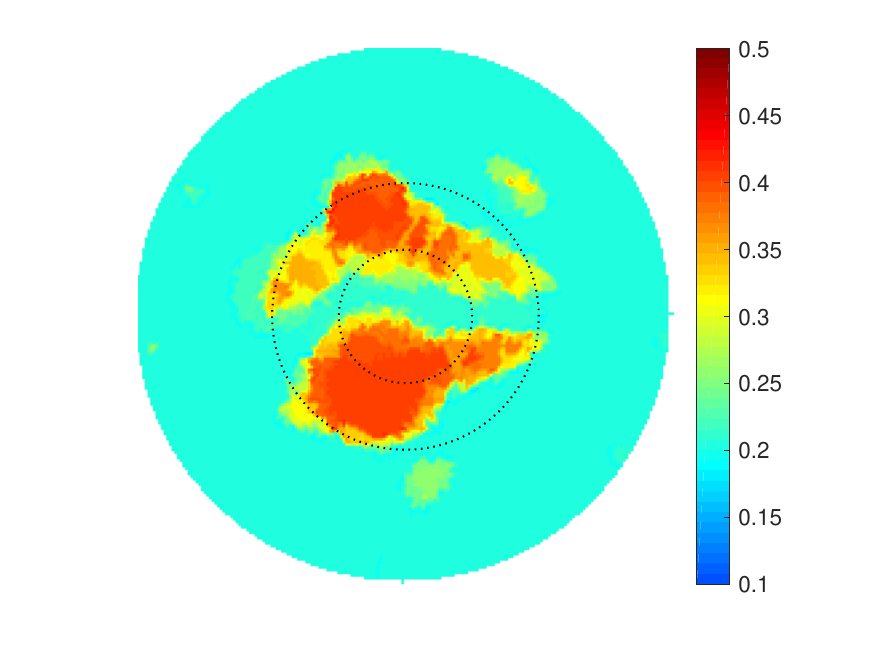}}}
  \end{center}
  \caption{Model~\#3. (a)~True electrical conductivity $\sigma_{true}(x)$.
    (b-f)~Solution images obtained by (b,c,e,f)~the
    proposed framework and (d)~the gradient-based SNOPT with PCA with
    (b,d,e)~no noise added and (c,f)~0.5\% noise in measurements. Images
    in (b,c,e,f) are obtained by utilizing (b,c)~$N_c^i = N_{c,max} = 8$
    and (e,f)~$N_c^i = N_{c,max} = 20$ conditions for the number of circles.}
  \label{fig:model_15}
\end{figure*}

Our last model~\#3, the hardest one, is created to check the method's
performance when the reconstructed region is not of a circular shape,
see Figure~\ref{fig:model_15}(a) with a C-shape region. As seen in
Figure~\ref{fig:model_15}(d), the gradient-based method with PCA is
unable to get a clear image even without noise, as the PCA
transform honors the structures of samples in the
$\mC(10000)$ collection. On the other hand, our new framework could
provide a good quality image with no noise in the data. This
performance may be further enhanced, even in the presence of noise,
once we re-set the number of circles condition to
$N_c^i = N_{c,max} = 20$. This proves the potential of the proposed
algorithm in applications with rather complex models.

\section{Concluding Remarks}
\label{sec:remarks}

In this work, we presented a novel computational approach for optimal
reconstruction of binary-type images useful in various applications
for (bio)medical practices. The proposed computational framework
uses a derivative-free optimization algorithm supported by a set
of sample solutions generated synthetically based on prior knowledge
of the simulated phenomena. This framework has an easy-to-follow design
tuned by a nominal number of computational parameters. High computational
efficiency is achieved by applying the coordinate descent method
customized to work with individual controls in the predefined custom order.

We investigated the performance of the complete framework in applications
to the 2D IPCD by the EIT technique. We claim, based upon our results,
that the proposed methodology of determining whether certain material
or medium contains a defective region has superior efficiency in comparison
with commonly used gradient-based and derivative-free techniques utilizing
control space parameterization via PCA. It is due primarily to the predominately
geometric nature of our approach, wherein we perturb known solutions
to similar related problems in order to converge to the best available
local/global minima.

There are many ways in which our proposed optimization framework can
be further tested and extended. Among other directions, we see the greatest
importance to test this approach in various applications to real data and
different types of defects/abnormalities in the media, as this would
certainly suggest new areas in which future developments may be required.
Even though this computational approach is initially tested
with synthetic EIT-related problems, we believe that this methodology could
also be easily applied to a broad range of problems seen in various fields
such as physics, geology, chemistry, etc.

\section*{Acknowledgements}
We wish to thank the anonymous reviewers for their valuable comments and
suggestions to improve the clarity of the presented approach and the overall
readability of this paper.

\bibliographystyle{spmpsci}
\bibliography{biblio_Bukshtynov,biblio_EIT,biblio_OPT}

\end{document}